\documentclass[preprint,12pt]{ourarticle}

\usepackage{graphicx}
\usepackage{natbib}
\usepackage{amsmath,amssymb,eucal}
\usepackage{mathrsfs}
\usepackage{stmaryrd}
\usepackage{color}
\usepackage{gensymb}
\usepackage[T1]{fontenc}
\usepackage{moreverb}
\usepackage{scalerel}
\usepackage{times}

\DeclareMathAlphabet{\mathpzc}{OT1}{pzc}{m}{it}
\usepackage{caption}
\newcommand{\eps}{\varepsilon}
\newcommand{\sig}{\sigma}

\newcommand{\tou}[1]{{\boldsymbol{#1}}}
\newcommand{\tod}[1]{{\boldsymbol{#1}}}
\newcommand{\toq}[1]{{\boldsymbol{#1}}}

\newcommand{\abs}[1]{\left|{#1}\right|}
\newcommand{\norm}[1]{\left|\left|{#1}\right|\right|}
\newcommand{\eff}[1]{\widetilde{#1}}

\newcommand{\ajj}[1]{\textcolor{black}{#1}}

\newcommand{\moy}[1]{\langle #1 \rangle}

\newcommand{\G}{\toq{\Gamma}}
\newcommand{\Gc}{\toq{\Gamma}_{\chi^{(1)}}}
\newcommand{\Gcr}{\toq{\Gamma}_{\chi^{(r)}}}
\newcommand{\Gcp}{\toq{\Gamma}_{\chi'}}
\newcommand{\Guc}{\toq{\Gamma}_{\chi^{(2)}}}
\newcommand{\Hh}{\toq{H}}
\newcommand{\Hc}{\toq{H}_{\chi^{(1)}}}

\newcommand{\Hcp}{\toq{H}_{\chi'}}
\newcommand{\Huc}{\toq{H}_{\chi^{(2)}}}
\newcommand{\mmoy}[1]{\langle{\langle #1 \rangle}\rangle}

\def\doublelow#1{\,\vtop{\ialign{\hfil$##$\hfil\crcr
                 \mathstrut #1 \crcr}}\,}

\newcommand{\h}{\toq{\mathcal{H}}}
\newcommand{\w}{\toq{\mathcal{W}}}
\date {}
\begin{document}

\begin{frontmatter}

\title{Convergence of iterative methods based on Neumann series for composite materials: theory and practice.
}
\author[label1]{Herv\'e Moulinec } 
\author[label1]{Pierre Suquet\corref{cor1}}
\author[label3]{Graeme W. Milton}
\address[label1]{Aix-Marseille Univ,  CNRS,  Centrale Marseille, LMA, \\ 4 Impasse Nikola Tesla, CS 40006,  13453 Marseille Cedex 13, France. }
\address[label3]{
    Univ Utah, Dept Math, Salt Lake City, UT 84112, USA.} 
\cortext[cor1]{Corresponding author: suquet@lma.cnrs-mrs.fr}

\begin{abstract}
Iterative Fast Fourier Transform methods are useful for calculating the fields in composite materials and their macroscopic response. By iterating back and forth until convergence,  the differential constraints are satisfied in Fourier space, and the constitutive law in real space.  The methods correspond to 
series expansions of appropriate operators and to series expansions for the effective tensor as a function of the component moduli. It is shown that the singularity structure of this function can shed much light on the convergence properties of the iterative Fast Fourier Transform methods. We look at a model example of a square array of conducting square inclusions  for which there is an exact formula for the effective conductivity (Obnosov). Theoretically some of the methods converge when the inclusions have zero or even negative conductivity. However,  the numerics do not always confirm this extended range of convergence and show that accuracy is lost
after relatively few iterations. There is little point in iterating beyond this. Accuracy improves
when the grid size is reduced, showing that the discrepancy is linked to the discretization. Finally, it is shown that none of the three iterative schemes investigated  over-performs the others for all possible microstructures and all contrasts.
 
\end{abstract}

  \end{frontmatter}

\section{Introduction}
It is a well-established result in the theory of linear composites that the local fields in such materials satisfy an integral equation  which can be  written symbolically as (Kröner \cite{KRO72}, Willis \cite{WIL81}, Milton \cite{MIL02})
\begin{equation}
  \left(\toq{I} + \toq{\Gamma}^0 \toq{\delta L}\right) \tod{\eps} = \tod{E},
  \label{LS}
\end{equation}
where 
$\tod{\eps}$ is the local field under investigation,  $ \tod{E}$ is its average, $\toq{\Gamma}^0$ is the  Green's operator associated with a homogeneous reference medium $\toq{L}^0$ and $\toq{\delta L}= \toq{L}-\toq{L}^0$ is the deviation of the actual material properties of the composite from the  homogeneous reference medium. As recognized by Kröner \cite{KRO72}, this integral equation has the same form as the Lippmann-Schwinger equation of scattering theory.
\vskip 0.1cm
The resolution of {(\ref{LS})} requires the inverse of $\toq{I} + \toq{\Gamma}^0 \toq{\delta L}$ which can be expanded in Neumann-Liouville series. This expansion was used by Brown \cite{BRO55} with one of the phases as reference material (see also Kröner \cite{KRO72}). 
\vskip 0.1cm
The Neumann-Liouville series is again the basis of the fixed-point scheme proposed by Moulinec and Suquet \cite{MOU94,MOU98} which differs from previous works by two important aspects. First, the reference medium is not chosen to be one of the phases, but optimized  to enhance the convergence of the series. Second, the {Discrete Fourier Transform} is used to compute efficiently the Green's operator  $\toq{\Gamma}^0$ and solve iteratively the Lippmann-Schwinger equation (\ref{LS}). This fixed-point algorithm (sometimes called the ``basic'' algorithm) alternates between real space (where the unit-cell is discretized along a regular grid) and Fourier space (where the spatial frequency vector takes discrete values dual to the spatial discretization in real space). The Fourier transform of the Green's operator $\toq{\Gamma}^0$ is known for  $\toq{L}^0$ with arbitrary anisotropy  and can be given {an} explicit form for several classes of symmetry (Khatchaturyan \cite{KHA83}, Mura \cite{MUR87}, Nemat-Nasser {\em et al} \cite{NEM82}). The choice of the reference tensor $\toq{L}^0$ and the rates of convergence that were observed with this fixed-point algorithm for various microstructures were consistent with the conditions for ``unconditional'' convergence derived by Michel {\em et al} \cite{MIC01}. 
\vskip 0.1cm
However, several issues were raised subsequently by the authors themselves or by other authors concerning the rate of convergence of the basic scheme, which can be poor when the contrast between the phases becomes large. In addition, spurious oscillations of the local fields were sometimes observed. The criticisms, and some of the progress made over the years, can be schematically listed under three main categories.  
\begin{enumerate}
\item
Rate of convergence. To improve on the rate of convergence of the basic scheme, accelerated methods have been proposed. The earliest one is due to Eyre and Milton \cite{EYR99} (generalized in Milton \cite{MIL02}) and consists in a change of the field on which the algorithm operates, {and a shift of $\toq{\Gamma}^0$ in an effort to reduce the norm of the operator product
entering the interations.} This accelerated scheme can also be seen as the summation of a Neumann-Liouville series. {The shifting and associated series expansion were 
first introduced in \cite{Milton:1990:RCF}, sect.5, {in the context of the conductivity problem, and its convergence properties, without assuming the conductivity tensor was symmetric,
were studied in \cite{Clark:1994:MEC}, sect.3.}
 The series is related to a series that forms the basis for establishing microstructure
independent relations satisfied by effective tensors (see, e.g., sections 14.10 and 17.3 in \cite{MIL02} and \cite{Grabovsky:2016:CMM})}.
Let us mention  that several other accelerated schemes (polarization scheme {\cite{MON12b,TO17}}, conjugate gradient \cite{ZEM10}, Galerkin approaches \cite{BRI12}) have also been proposed but are not discussed here. {Some comparisons were recently made in \cite{Mishra:2016:CSL}, suggesting the conjugate gradient method may often have the fastest convergence.}
\item
Convergence criterion. As in every iterative method, choosing a sensible test to decide when the iterations should be stopped is of crucial importance. The initial criterion of Moulinec and Suquet \cite{MOU98} was based on the $L^2$ norm of one of the equations to be satisfied (equilibrium equation). Other criteria, such as the difference between two iterates, have also been proposed (Milton \cite{MIL02}).  
\item
Spectral derivatives and Green's operators. The Green's operator used in the basic scheme is the continuous Green's operator. Spurious oscillations in the solutions have been observed in the computed fields. Gibbs phenomena are sometimes invoked to explain these oscillations, but in our opinion their origin is different, as they appear gradually along the iterations. A possible explanation is that computing the Fourier transform of the derivative of a nonsmooth function by means of the Fourier transform requires some care. This has motivated M\"{u}ller \cite{MUL96} and subsequently several authors ( Brown {\em et al} \cite{BRO02}, Willot \cite{WIL14}, Schneider {\em et al} \cite{SCH16} among others) to use discrete Green's operators.   
\end{enumerate}
The present article discusses several aspects of the convergence of three Neumann-Liouville series for three different schemes available in the literature and presented here in a unified form. \ajj{How power series enter algorithms for determining the effective properties of composite materials is recalled in section \ref{Neumann},  together with theoretical results on the convergence of the series with no restriction on the microstructure. These series expansions and the corresponding sufficient conditions for convergence are valid for general $N$-phase composites.} The integral operator entering the power series is split into an integral operator with norm 1 and a local operator depending only on the contrast between the constituents. Section 3 examines the theoretical convergence of the series when more information about the microstructure is available. The integral operator is {split} differently into an integral operator depending on the microstructure and a local operator depending only on the contrast between the phases. Theoretical estimates for the radius of convergence of the iterative methods are derived. \ajj{Again these results are valid for $N$-phase composites.} \ajj{For two-phase composites, where the contrast between the phases can be measured with a single parameter $z$}, it is shown that information on the singularities of the effective moduli in the complex $z$-plane can be used to improve the range of phase contrast for which these iterative methods converge and the theoretical prediction of their rate of convergence. These theoretical estimates are compared in section 4 with numerical simulations of  the effective conductivity of a square array of square inclusions for which an explicit solution is available (Obnosov \cite{OBN99}). {The range of contrast for which convergence is observed numerically appears to be smaller than what it should be theoretically. It is found that the origin of this discrepancy is the reduction of the continuous problem to a finite-dimensional one. The discretized solution being different from the exact continuous one, the higher- order terms in the series have to be different from their exact value. Finally the rate of convergence of the three iterative schemes is examined for more general microstructures when the location of the singularities and the contrast between the phases are varied. It is found that none of the three iterative schemes over-performs the other two for all possible microstructures and all phase contrasts}. 

\section{Estimates for the convergence radius of Neumann-Liouville series with no information on the microstructure}
\subsection{The Lippmann-Schwinger equation}
Consider a unit-cell $V$ of a periodic composite material. The composite is made up of  $N$ homogeneous phases $V^{(r)}, r=1,..,N$, whose distribution is defined by characteristic functions 
$\chi^{(r)}$. {L}et $\moy{.}$ and $\moy{.}^{(r)}$ denote spatial averaging over $V$ and $V^{(r)}$ respectively. The material property under consideration is characterized by a tensor field $\toq{L}(\tou{x})$ relating two local fields $\tod{\sig}(\tou{x})$ and $\tod{\eps}(\tou{x})$ satisfying partial differential equations expressing, in elasticity, balance equations and compatibility conditions :
\begin{equation}
	\left. 
        \begin{array}{l}
\tod{\sigma}(\tod{x})=\displaystyle\toq{L}(\tou{x}):\tod{\eps}(\tou{x}), \quad
 \hbox{\rm div }\tod{\sigma}(\tou{x})=0 \quad \hbox{\rm in }\ V , \\[2ex]
 \tod{\eps}(\tou{x}) = {\tod{E}} + \tod{\eps}^*(\tou{x}), \quad  \tod{\eps}^*(\tou{x}) = 
\frac{1}{2} ( \tod{\nabla} \tou{u}^*(\tou{x}) + \tod{\nabla} {\tou{u}^*}^T(\tou{x})),  \\[2ex] 
 \tou{u}^* \ \text{periodic on}\ \partial V, \quad \tod{\sig}.\tou{n} \ \text{anti-periodic on}\ \partial V.
 \end{array}
 \right\}
 \label{local1}
 \end{equation}
In elasticity $\toq{L}$ is the fourth-order stiffness tensor, $\tod{\sig}$ is the stress field, $\tod{\eps}$ is the strain field, while in conductivity, $\toq{L}$ is the second-order conductivity tensor, $\tod{\sig}$ is the current field, $\tod{\eps}$ is the electric field which is the gradient of the electrical potential $u$. $\toq{L}$ is assumed to be a symmetric tensor.
\vskip 0.1cm
Introducing a reference tensor $\toq{L}^0$ (symmetric and definite positive), the constitutive relation between $\tod{\sig}$ and $\tod{\eps}$ is {rewritten} with a polarization field $\tod{\tau}$ as 
\begin{equation}
\tod{\sigma}(\tod{x})=\displaystyle \toq{L}^0:\tod{\eps}(\tou{x}) + \tod{\tau}(\tou{x}), \quad \tod{\tau}(\tou{x})=\toq{\delta L}(\tou{x}):\tod{\eps}(\tou{x}), \quad   \toq{\delta L}(\tou{x})= \toq{L}(\tou{x})-\toq{L}^0.
\label{tau}
\end{equation}
The solution of (\ref{local1}) can be expressed with the Green's operator associated with $\toq{L}^0$ as 
$$  \tod{\eps} = \tod{E} - \toq{\Gamma}^0 \tod{\tau}, $$
where $\toq{\Gamma}^0$ is an integral operator  {(see appendix A)} and $\toq{\Gamma}^0 \tod{\tau}$ is to be understood as the action of this integral operator on the field $\tod{\tau}$ given by (\ref{tau}). Replacing $\tod{\tau}$ by its expression (\ref{tau}) leads to the Lippmann-Schwinger integral equation  {for the field $\tod{\eps}$}:
\begin{equation}
\left(\toq{I} + \toq{\Gamma}^0 \toq{\delta L}\right) \tod{\eps} = \tod{E}, 
\label{LS1}
\end{equation}
{the solution of which can be written as 
\begin{equation}
 \tod{\eps} = \left(\toq{I} + \toq{\Gamma}^0 \toq{\delta L}\right)^{-1}\tod{E}. 
\label{LS2a}
\end{equation}}
 {The effective moduli $\eff{\toq{L}}$ relating the average stress $\moy{\tod{\sig}}$ and the average strain $\tod{E}=\moy{\tod{\eps}}$ read as
\begin{equation}
\eff{\toq{L}} = \moy{\toq{L} \left(\toq{I} + \toq{\Gamma}^0 \toq{\delta L}\right)^{-1} } = \toq{L}^0 + \moy{  \toq{\delta L} \left(\toq{I} + \toq{\Gamma}^0 \toq{\delta L}\right)^{-1}}.
\label{LS2b}
\end{equation}
}
\subsection{Neumann-Liouville series}\label{Neumann}
The operator $(\toq{I} + \toq{\Gamma}^0 \toq{\delta L})^{-1}$ can be {\em formally} expanded in power series of $\toq{\Gamma}^0 \toq{\delta L}$ {and the corresponding expansion for the strain field $\tod{\eps}$ and the effective moduli $\eff{\toq{L}}$ are:}
\begin{equation}
 \tod{\eps} = \sum_{j=0}^\infty \left(- \toq{\Gamma}^0 \toq{\delta L}\right)^{j}\tod{E}, \quad 
 \eff{\toq{L}}= \toq{L}^0 + \sum_{j=0}^\infty \moy{  \toq{\delta L} \left(-\toq{\Gamma}^0 \toq{\delta L}\right)^{j}}.
\label{LS3}
\end{equation}
The successive terms in these series correspond to successive iterates in the resolution of the Lippmann-Schwinger equation by Picard iterations \footnote{Here (as in Milton \cite{MIL02} chap. 14) $\toq{\Gamma}^0$ and $\toq{\delta L}$ should be interpreted as {\em operators} and  $\left(- \toq{\Gamma}^0 \toq{\delta L}\right)^{j}$ as the operator  $- \toq{\Gamma}^0 \toq{\delta L}$ applied $j$ times (and not as the field $- \toq{\Gamma}^0 \toq{\delta L}$ raised to the $j$-th power).}
\begin{equation}
 \tod{\eps}^{(k+1)} = - \toq{\Gamma}^0 \toq{\delta L}\tod{\eps}^{(k)} + \tod{E}, \; \text{where} \;  \tod{\eps}^{(k)} = \sum_{j=0}^k \left(- \toq{\Gamma}^0 \toq{\delta L}\right)^{j}\tod{E} .
\label{LS4}
\end{equation}
Similarly the $k$-th approximation of the effective tensor $\eff{\toq{L}}$ is given as 
\begin{equation}
 \eff{\toq{L}}^{(k)}= \toq{L}^0 + \sum_{j=0}^k \moy{  \toq{\delta L} \left(-\toq{\Gamma}^0 \toq{\delta L}\right)^{j}}. 
\label{LS5}
\end{equation}
Convergence of the Liouville-Neumann series (\ref{LS3}), or equivalently of the iterative procedure (\ref{LS4}), requires specific choices for the reference tensor $\toq{L}^0$ and/or limitations on the contrast ${\toq{L}^0}^{-1}\toq{\delta L}$ between the phases. To investigate these conditions it is useful to re-write $\toq{\Gamma}^0 \toq{\delta L}$ as the composition of two operators
{
\begin{equation}
\toq{\Gamma}^0 \toq{\delta L} = \toq{\Gamma}^0 \toq{L}^0 \ {\toq{L}^0}^{-1} \toq{\delta L} = \toq{\Gamma}^1  \toq{Z}, \quad \toq{\Gamma}^1=\toq{\Gamma}^0 {\toq{L}^0},
\quad \toq{Z} = {\toq{L}^0}^{-1} {\toq{\delta L}}.
\label{split}
\end{equation}
}
$\toq{\Gamma}^1$ is a projection (for an energetic scalar product, see appendix A) and its norm is therefore equal to 1. The power series (\ref{LS3})$_1$ and {(\ref{LS3})$_2$} take the form
\begin{equation}
 \tod{\eps} = \sum_{j=0}^\infty \left(- \toq{\Gamma}^1 \toq{Z}\right)^{j}\tod{E}, \quad { \toq{L}^0 }^{-1}  \eff{\toq{L}}= \toq{I} + \sum_{j=0}^\infty \moy{  \toq{Z} \left(-\toq{\Gamma}^1 \toq{Z}\right)^{j}},
 \label{LS3b}
\end{equation}
where it is seen that, in this approach, the two pivotal operators are $\toq{\Gamma}^1$ and ${\toq{Z}}$. 
\vskip 0.1cm
Eyre and Milton \cite{EYR99} noticed that the power series (\ref{LS3}) makes use of powers of the operator $\toq{\Gamma}^1=\toq{\Gamma}^0 \toq{L}^0$. Given that this operator is positive with norm equal to 1 (in energetic norm), they suggested to iterate with a shifted operator and introduced (with different notations)
\begin{equation}
\toq{H}^1 = 2 \toq{\Gamma}^0 \toq{L}^0 - \toq{I}.
\label{shift}
\end{equation}
$\toq{H}^1$ is  a self-adjoint operator (for an appropriate scalar product), with all its eigenvalues between $-1$ and $1$ and can therefore replace $\toq{\Gamma}^1$ in the Lippmann-Schwinger equation. Following Milton \cite{MIL02} an equivalent form of the Lippmann-Schwinger equation and of its associated Neumann-Liouville series can be derived after some algebra. Noting  that 
$\toq{\Gamma}^0 \toq{L}^0 = \frac{1}{2}(\toq{H}^1+ \toq{I})$, one obtains that 
\begin{equation}
 \toq{I}+\toq{\Gamma}^0 \toq{\delta L}= \toq{I}+ \frac{1}{2}(\toq{H}^1+ \toq{I}) {\toq{L}^0}^{-1} \toq{\delta L} =
  \toq{I}+ \frac{1}{2} {\toq{L}^0}^{-1} \toq{\delta L} + \frac{1}{2}\toq{H}^1 {\toq{L}^0}^{-1} \toq{\delta L}.
 \label{upsilon0}
 \end{equation}
  Note that 
  \begin{equation}
    \begin{array}{ll}
      \toq{I}+ \frac{1}{2} {\toq{L}^0}^{-1}\toq{\delta L}
      &= {\toq{L}^0}^{-1} {\toq{L}^0} + \frac{1}{2} {\toq{L}^0}^{-1}\toq{\delta L} \\
      &= \frac{1}{2} {\toq{L}^0}^{-1}(\toq{L}+ \toq{L}^0) , \\
    \end{array}
    \label{upsilon2a}
  \end{equation}
  and that 
  $$
    \begin{array}{ll}
      {\toq{L}^0}^{-1} \toq{\delta L}  \left(\toq{L}+\toq{L}^0\right)^{-1} \toq{L}^0
   & = {\toq{L}^0}^{-1} \left( \toq{L} + \toq{L}^0 - 2 \toq{L}^0 \right) \left(\toq{L}+\toq{L}^0\right)^{-1} \toq{L}^0 
      = \toq{I} - 2 \left(\toq{L}+\toq{L}^0\right)^{-1} \toq{L}^0 \\ & = \left(\toq{L}+\toq{L}^0\right)^{-1} \left(\toq{L}+\toq{L}^0\right) - 2 \left(\toq{L}+\toq{L}^0\right)^{-1} \toq{L}^0  =  \left(\toq{L}+\toq{L}^0\right)^{-1} \toq{\delta L}. 
    \end{array}
   $$
   and therefore 
   \begin{equation}
   {\toq{L}^0}^{-1} \toq{\delta L} = \toq{W}  {\toq{L}^0}^{-1} \left(\toq{L}+\toq{L}^0 \right) , \quad  
   \frac{1}{2}\toq{H}^1 {\toq{L}^0}^{-1} \toq{\delta L} =  \frac{1}{2}\toq{H}^1 \toq{W} \frac{1}{2} {\toq{L}^0}^{-1}(\toq{L}+ \toq{L}^0)\quad \text{with} \quad \toq{W} = 
 \left(\toq{L}+\toq{L}^0\right)^{-1} \toq{\delta L}.
     \label{upsilon2b}
  \end{equation}
Substituting {\eqref{upsilon2a} and \eqref{upsilon2b}}  into  \eqref{upsilon0} yields
\begin{equation} \toq{I}+\toq{\Gamma}^0 \toq{\delta L}= (\toq{I}  + \toq{H}^{1} \toq{W}) \frac{1}{2} {\toq{L}^0}^{-1}(\toq{L}+ \toq{L}^0). 
\label{upsilon3}
\end{equation}
The Lippmann-Schwinger equation (\ref{LS1}) and its associated Neumann series (\ref{LS3}) become
 \begin{equation}
  \left( \toq{I}+\toq{H}^1 \toq{W} \right)
  \frac{1}{2} {\toq{L}^0}^{-1}(\toq{L}+ \toq{L}^0) \tod{\eps} = \toq{E},   \quad \tod{\eps} = 2  (\toq{L}+ \toq{L}^0)^{-1} \toq{L}^0 \sum_{j=0}^{\infty} (- \toq{H}^1 \toq{W}  )^{j} \tod{E}.
 \label{LS7}
 \end{equation}
The corresponding series expansion for the effective moduli is obtained by writing that 
$$\tod{\sig} = \toq{L}^0:\tod{\eps} + \toq{\delta L}:\tod{\eps}, \quad \text{i.e.} \quad \moy{\tod{\sig}} = \toq{L}^0:\tod{E} + \moy{\toq{\delta L}:\tod{\eps}}$$
and using again {\eqref{upsilon2a} and \eqref{upsilon2b}}, one gets
\begin{equation}
{\toq{L}^0}^{-1} \eff{\toq{L}} = \toq{I}  +  2 \sum_{j=0}^\infty \moy{\toq{W} \left(-\toq{H}^1  \toq{W}\right)^{j}},
\label{LS9}
\end{equation}
which evidences the similarity with (\ref{LS3b}), $\toq{Z}$
being substituted with 
$\toq{W}$, and  $\toq{\Gamma}^1$ with  $\toq{H}^1$.
\vskip 0.1cm
\noindent{\em Remark 1:}
Truncated sums corresponding to $\tod{\eps}^{(k)}$ in (\ref{LS4}) can also be considered within the Eyre-Milton scheme. Defining
\footnote{The notation $\toq{e}^{(k)}$ is used here instead of $\tod{\eps}^{(k)}$ to highlight the fact that the truncated sums are {not guaranteed} compatible strain fields (except at convergence), unlike the truncated series (\ref{LS4}) and (\ref{EMconf}). }
\begin{equation}
  \tod{e}^{(k)} =  2  (\toq{L}+ \toq{L}^0)^{-1} \toq{L}^0 \sum_{j=0}^{k} (- \toq{H}^1 \toq{W}  )^{j} \tod{E}, 
 \label{partial1a}
 \end{equation}
 and using \eqref{upsilon3}, the following recurrence relation is obtained
 $$ \left(\toq{I} + \toq{\Gamma}^0 \toq{\delta L} \right)\tod{e}^{(k)} = \left(\toq{I} + \toq{H}^1 \toq{W} \right) \sum_{j=0}^{k} (- \toq{H}^1 \toq{W}  )^{j} = 
\toq{E}+  \frac{1}{2} {\toq{L}^0}^{-1}(\toq{L}+ \toq{L}^0) \left( \tod{e}^{(k)} -  \tod{e}^{(k+1)}\right),
 $$ 
 or equivalently,
\begin{equation}
 \tod{e}^{(k+1)} =   \tod{e}^{(k)} - 2   (\toq{L} + \toq{L}^0)^{-1} \toq{L}^0 \left[ (\toq{I} + \toq{\Gamma}^0 \toq{\delta L} )\tod{e}^{(k)} - \tod{E}\right],
 \label{partial2}
 \end{equation}
 which is equivalent to the writing of the Eyre-Milton algorithm in Michel {\em et al}  \cite{MIC01}.
\vskip 0.1cm
{\noindent{\em Remark 2:}
Noting that $\toq{\Gamma}^0\toq{\delta L}\tod{\eps}=\tod{E}-\tod{\eps}$, the last equation in (\ref{LS7}) can be rewritten as
\begin{equation}
\tod{\eps}=\tod{E}-2\toq{\Gamma}^0\toq{W}\toq{L}^0 \sum_{j=0}^{\infty} (- \toq{H}^1 \toq{W}  )^{j} \tod{E}.
\end{equation}
The partial sums 
\begin{equation}
 \label{EMconf}
\tod{\eps}^{(k)}=\tod{E}-2\toq{\Gamma}^0\toq{W}\toq{L}^0 \sum_{j=0}^{k} (- \toq{H}^1 \toq{W}  )^{j} \tod{E}
=\tod{E}-\toq{\Gamma}^0\toq{\delta L}\tod{e}^{(k)},
\end{equation}
then define a sequence of conforming approximations to $\tod{\eps}$, in the sense that each field $\tod{\eps}^{(k)}$ is a compatible strain field, i.e.,
the gradient of some symmetrized displacement gradient.
Note that in the process of computing the iterates in (\ref{partial2}) one calculates $\toq{\Gamma}^0 \toq{\delta L}\tod{e}^{(k)}$ thus 
immediately giving $\tod{\eps}^{(k)}$ from (\ref{EMconf}) with no extra work. In conclusion, it is inaccurate to say that the Eyre-Milton algorithm is a non-conforming approximation scheme
\textemdash one only needs to keep track of the fields $\tod{\eps}^{(k)}$.}
\vskip 0.1cm
\noindent{\em Remark 3:} An alternative Neumann series for the Eyre-Milton scheme can be obtained {by writing the iterative algorithm (\ref{partial2})} in terms of the polarization field 
\begin{equation}
\tod{\tau} = (\toq{L}+\toq{L}^0) \tod{e}.
\label{polar}
\end{equation}
 Relation (\ref{partial2}), re-written as 
\begin{equation}
  ( \toq{L} + \toq{L}^0 ) {\tod{e}}^{k+1}
  =
  ( \toq{L} + \toq{L}^0 ) {\tod{e}}^{k}
  - 2 \toq{L}^0 {\tod{e}}^{k}
  - 2 \toq{L}^0 \toq{\Gamma}^0
  \left( \toq{L} - \toq{L}^0 \right) {\tod{e}}^{k} 
  + 2 \toq{L}^0 \tod{E} \ ,
\label{MBscheme3}
\end{equation}
can be expressed with the polarization field as 
\begin{equation}
    {\tod{\tau}}^{k+1}
    =
    - \h^1
    \w {\tod{\tau}}^{k} 
    + 2 \toq{L}^0 \tod{E} \ ,\;  \h^1 =   2 \toq{L}^0 \toq{\Gamma}^0 -   \toq{I}, \;    \w = ( \toq{L} - \toq{L}^0) (\toq{L}+\toq{L}^0)^{-1}.
\label{EM3}
\end{equation}
Note that the operators $\h^1$ and $\w$  are slightly different from the operators  $\toq{H}^1$ {and} $\toq{W}${: see}
 (\ref{shift}) and (\ref{upsilon3}).
\vskip 0.1cm
The series expansion corresponding to (\ref{EM3}) now reads 
\begin{equation}
  {\tod{\tau}}^{k}
  =
  \sum_{k=0}^{k}
2  \left( - \h^1 \w \right)^k \  \toq{L}^0 \tod{E} \ .
\label{EM4}
\end{equation}
At convergence $\tod{e}=\tod{\eps}$ and the following two expressions for $\moy{\tod{\tau}}$ resulting from (\ref{polar}) and (\ref{EM4}){,}
\begin{equation}
  {\tod{T}} = < \tod{\tau} > =
  \sum_{k=0}^{\infty}
2  < \big( - \h^1 \w )^k > \ \toq{L}^0 \tod{E}
  \nonumber
\end{equation}
and 
\begin{equation}
  \tod{T} = < \tod{\tau} > = \toq{L}^0 <\tod{\eps}> + <\tod{\sigma}>
  = \toq{L}^0 \tod{E} + \tilde{\toq{L}} \tod{E}{,}
\end{equation}
lead to
\begin{equation}
    \tilde{\toq{L}} { \toq{L}^0}^{-1} = 
    \toq{I} +  \sum_{k=1}^{\infty}
    2 \moy{ \big( - \h^1 \w )^k }. 
 \label{Tau_eff}
\end{equation}
{Noting that 
\begin{equation}
 \moy{\big( \h^1 \w )^k } =  \moy{ (2 \toq{L}^0 \toq{\Gamma}^0 -   \toq{I}) \w \big( \h^1 \w )^{k-1} }=
  2 \toq{L}^0  \moy{\toq{\Gamma}^0 \w \big( \h^1 \w )^{k-1} }- \moy{\w \big( \h^1 \w )^{k-1} }= 
 - \moy{\w \big( \h^1 \w )^{k-1} }, 
 \label{Tau_eff2}
\end{equation}
 one gets from (\ref{Tau_eff})
\begin{equation}
    \tilde{\toq{L}} { \toq{L}^0}^{-1} = 
    \toq{I} +  2  \sum_{j=0}^{\infty}
     \moy{ \w \big( - \h^1 \w )^{j} }, 
 \label{Tau_eff3}
\end{equation}
which is similar to (\ref{LS9}).}

\subsection{Unconditional  convergence of the power series}\label{unconditional}
"Unconditional" refers to sufficient conditions which are  independent of the microstructure under consideration. Conservative estimates on the range of contrast for which the power series (\ref{LS3})$_1$ converges,  can be obtained by estimating the norm (defined in an appropriate way) of $\toq{\Gamma}^0 \toq{\delta L}$. This can be done in a first attempt by using the decomposition (\ref{split}). 
\vskip 0.1cm
It has been already recalled that the operator $\toq{\Gamma}^1$ has norm equal to 1  independently of the microstructure (appendix A). Therefore 
\begin{equation}
\norm{\toq{\Gamma}^0 \toq{\delta L}} \leq \norm{ {\toq{L}^0}^{-1} \toq{\delta L}}.
\label{CV1}
\end{equation}
A sufficient condition for convergence of the power series is therefore that 
\begin{equation}
\norm{ {\toq{L}^0}^{-1} \toq{\delta L}} < 1 .
\label{CV2}
\end{equation}
When the reference medium is chosen to be one of the phases, as proposed by Brown \cite{BRO55}, the inequality (\ref{CV2}) imposes a restriction on the moduli of the other phases which, roughly speaking, have to be
{lower than twice}
the moduli of the phase taken as reference. The corresponding scheme (\ref{LS3}) will be called B-scheme (B for Brown) in the sequel.
\vskip 0.1cm
Moulinec and Suquet \cite{MOU94} suggested to consider $\toq{L}^0$ as a free parameter and to choose it in a such a way that (\ref{CV2}) is satisfied. There is a wide range of possible reference media. For instance taking $\toq{L}^0$ greater (in the sense of quadratic forms) than all the individual moduli $\toq{L}^{(r)}, r=1,..,N$ of the phases ensures that (\ref{CV2}) is satisfied with no additional restriction on the contrast between the phases.
However, even when the convergence of the power series (\ref{LS3}) is ensured, its convergence rate, which  is actually strongly dependent on  the reference medium, can be low. This convergence rate is directly related to the spectral radius of $\toq{\Gamma}^0 \toq{\delta L}$ which is bounded from above by  $\norm{ {\toq{L}^0}^{-1} \toq{\delta L}}$ independently of the microstructure. The ``unconditionally'' optimal reference medium (with no information on the microstructure) is obtained by solving the optimization problem 
\begin{equation}
\doublelow{\text{Min} \cr
\toq{L}^0 } \norm{{\toq{L}^0}^{-1}:\toq{\delta L}},
\label{CV2b}
\end{equation}
or equivalently for a $N$-phase composite, the discrete optimization problem 
\begin{equation}
\doublelow{\text{Min} \cr
\toq{L}^0 }  \doublelow{\text{Max} \cr
r=1,...,N }\norm{{\toq{L}^0}^{-1}:(\toq{L}^{(r)}- \toq{L}^0)}.
\label{CV3}
\end{equation}
This optimization problem leads for instance in a two phase composite with isotropic phases to  (Moulinec and Suquet \cite{MOU98}, Michel {\em et al} \cite{MIC01}, Milton \cite{MIL02})
\begin{equation}
\toq{L}^{(0)} = \frac{1}{2} (\toq{L}^{(1)} + \toq{L}^{(2)}).
\label{CV4}
\end{equation}
The corresponding scheme (\ref{LS3}) with the choice (\ref{CV4}) will be referred to as the MS-scheme {(for Moulinec-Suquet)} in the sequel.
\vskip 0.1cm
Similar sufficient conditions can be obtained for the Eyre-Milton version (\ref{LS7}) of the Lippmann-Schwinger equation. The operator $\toq{H}^1$ has norm 1 (with an appropriate choice of the scalar product), and the spectral radius of $\toq{H}^1 \toq{W}$ satisfies
 \begin{equation}
\norm{\toq{H}^1 \toq{W}} \leq \norm{\toq{W}}.
\label{CV1b}
\end{equation}
Therefore convergence of the power series (\ref{LS7})$_2$ is assured whenever
 \begin{equation}
\norm{\toq{W}}= \norm{(\toq{L}+ \toq{L}^0)^{-1}(\toq{L}- \toq{L}^0)}< 1.
\label{condconv2}
\end{equation}
This condition is always satisfied when $\toq{L}$ and $\toq{L}^0$ are definite positive, bounded, and bounded from below by some positive constant of coercivity. Again a large choice of reference media is possible. The ``unconditionally'' optimal choice for $\toq{L}^0$ for the requirement (\ref{condconv2}) is obtained by minimizing the left-hand side (\ref{condconv2})
  \begin{equation}
\doublelow{\text{Min} \cr
\toq{L}^0 } \norm{(\toq{L}+ \toq{L}^0)^{-1}(\toq{L}- \toq{L}^0)}. 
\label{condconv2bis}
  \end{equation}
An explicit solution to this minimization problem can be derived when the tensors $\toq{L}$ and $\toq{L}^0$  have a more specific form, as will be illustrated  in section  \ref{simpl}. 
The  scheme (\ref{LS7}) with this optimized reference medium will be referred to as the {EM-scheme} (for Eyre-Milton) in the sequel. 
\vskip 0.1cm
The microstructure does not enter the optimization problems {(\ref{CV2b})} and (\ref{condconv2bis}). Solving these minimization problems ensures that the spectral radius of the two operators $\toq{\Gamma}^0\toq{\delta L}= \toq{\Gamma}^1\toq{Z}$  and $\toq{H}^1 \toq{W}$ is always smaller than 1, for any microstructure and any strictly positive and finite moduli in the phases. \vskip 0.1cm
However in the extreme case of one  phase being void (with zero moduli), the upper bounds (\ref{CV1}) and  (\ref{CV1b}) on the spectral radii of the iterative methods are exactly 1. Convergence is not ensured. The situation is even worse when the modulus  of one of the phases is nonpositive while the moduli of the other phases are positive. In the other extreme case, one of the phases having infinite moduli, the ``optimal'' reference medium itself has infinite moduli (and it does not make sense to consider its Green's operator $\toq{\Gamma}^0$).  
\vskip 0.1cm
A sharper estimate of this spectral radius can be obtained by taking into account the microstructure of the composite. This is the aim of section \ref{sec3} where the class of composites under consideration will be restricted in order to avoid complicated notations.
\subsection{A class of composite materials with a scalar measure of the contrast}\label{simpl}
In general, the material properties of the phases depend on several moduli, either because the phases are anisotropic or because the tensors describing the material properties under consideration, even when they are isotropic, depend on several moduli (isotropic elastic tensors depend on two moduli, a bulk modulus and a shear modulus). Therefore the contrast between the different phases depend on several contrast variables. 
\vskip 0.1cm 
The results of the forthcoming sections are better understood when the contrast between the different phases is measured with a single scalar parameter, called $z$ in the sequel. This is in particular the case when all tensorial moduli $\toq{L}^{(r)}$ are proportional to the same tensor ${\toq{\Lambda}}$ and when the moduli of the reference medium are also chosen proportional to the same tensor ${\toq{\Lambda}}$. For simplicity again, it will be sufficient for our purpose here to restrict attention to two-phase composites
\begin{equation}
\toq{L}^{(1)}= z {\toq{\Lambda}}, \quad \toq{L}^{(2)}=  {\toq{\Lambda}}, \quad \toq{L}^{0}= z_0 {\toq{\Lambda}}, 
\label{conduc1}
\end{equation}
where after a proper rescaling, the prefactor of ${\toq{\Lambda}}$ in the phases can be set equal to  1 in phase 2, to $z$ in phase 1 and to $z_0$ in the reference medium. 

\ajj{This is typically the case when investigating  the effective conductivity of two-phase composites with isotropic phases. Then
\begin{equation}
{\toq{\Lambda}}=  {\toq{I}},
\label{cas1}
\end{equation}
and $z$ is the conductivity of phase 1 (the conductivity of phase 2 being 1). This is the situation considered in section \ref{section_theoretical}. A similar situation is encountered when investigating the effective elastic moduli of two-phase composites with isotropic phases sharing the same Poisson ratio $\nu$, ${\toq{\Lambda}}= \frac{1}{1- 2\nu} \toq{J} + \frac{1}{1+\nu} \toq{K}$ where $\toq{J}$ and $\toq{K}$ are the usual projectors on second-order spherical and deviatoric symmetric tensors respectively, $z$ is the Young modulus of phase 1 (phase 2 having a normalized Young modulus equal to 1). The case of incompressible elastic phases considered in \cite{HOA14} is a special case of (\ref{conduc1}) with, as usual, some care to be taken to handle properly incompressible phases ($\nu=1/2$). }
\ajj{In two-dimensions (plane-strain or plane-stress) the assumption of equal Poisson ratios is in fact not restrictive at all since a problem with unequal
Poisson ratios and unequal Young moduli can be mapped using the ``CLM Theorem'' (see, for example, \cite{CLM} and equation (4.23) and (4.24) in \cite{MIL02}) 
to a mathematically equivalent elasticity problem where the Poisson ratios are equal (in this equivalent problem the elasticity tensors are not necessarily positive definite,
but this is not a problem as existence and uniqueness of solutions is guaranteed by the equivalence to the original problem where these tensors are positive definite). Also
we note that that the analysis for two-dimensional conductivity problems immediately extends to the equivalent problem of antiplane elasticity.
}
\vskip 0.1cm
Under {assumption (\ref{conduc1})} the operators $\toq{\Gamma}^1$, $\toq{H}^1$ are independent of the modulus $z_0$  of the reference medium, and they take a simple form, 
\begin{equation}
\toq{\Gamma}^1 = {\toq{\Gamma^{(\Lambda)}}} {\toq{\Lambda}}, \quad \toq{H}^1 = 2{\toq{\Gamma^{(\Lambda)}}} {\toq{\Lambda}} - \tod{I}, 
\label{simpl1}
\end{equation}
where  ${\toq{\Gamma^{(\Lambda)}}}$ is the Green's operator {$\toq{\Gamma}^0$} for  ${\toq{\Lambda}}$. The local tensor fields $\toq{Z}$ and $\toq{W}$ depend on the contrast $z$ as
\begin{equation}
\left.
\begin{array}{c}
\displaystyle
\toq{Z}(\tou{x}) = \left[\frac{z-z_0}{z_0}\chi^{(1)}(\tou{x}) + \frac{1-z_0}{z_0}\chi^{(2)}(\tou{x})\right] \toq{I},\\[2ex]
\displaystyle
\toq{W}(\tou{x}) =\left[ \frac{z-z_0}{z+z_0}\chi^{(1)}(\tou{x}) + \frac{1-z_0}{1+z_0}\chi^{(2)}(\tou{x})\right] \toq{I}.
\end{array}
\right\}
\label{simpl2}
\end{equation}
The operators $\toq{Z}$ and $\toq{W}$ are local  so their operator norm is a function of their pointwise norm in each phase
\begin{equation}
\norm{\toq{Z}} = \text{Max} \left(  \abs{\frac{z-z_0}{z_0}}, \abs{\frac{1-z_0}{z_0}} \right) \abs{\toq{I}}, \quad 
\norm{\toq{W}} = \text{Max} \left(  \abs{\frac{z-z_0}{z+z_0}}, \abs{\frac{1-z_0}{1+z_0}} \right) \abs{\toq{I}}, 
\label{simpl3}
\end{equation}
where $\abs{\toq{I}}$ is the Euclidian norm of the identity which can be set equal to $1$ by a proper rescaling. The convergence of three different series expansions will be investigated in the {remainder} of the paper. The first two series correspond to the expansion (\ref{LS3}) with the two variants mentioned above, B-scheme when  the reference medium is one of the phases, MS-scheme when  the reference medium  minimizes $\norm{\toq{Z}}$ (optimization problem (\ref{CV2b})). The third series is the EM-scheme obtained with the expansion (\ref{LS7}) and the reference medium minimizing $\norm{\toq{W}}$ (optimization problem (\ref{condconv2bis})). 
\begin{enumerate}
\item
{\em B-scheme:} When the reference medium is phase 2 ($z_0=1$), then  $\norm{\toq{Z}}= \abs{z-1}$ and the sufficient condition (\ref{CV2}) (independent of the microstructure) is 
\begin{equation}
\abs{z-1} < 1.
\label{simpl4}
\end{equation}
\item
{\em MS-scheme:} The optimal reference medium minimizing $\norm{\toq{Z}}$  is found by solving the minimization problem 
\begin{equation}
\doublelow{\text{Min} \cr
z_0} \ \text{Max} \left(  \abs{\frac{z-z_0}{z_0}}, \abs{\frac{1-z_0}{z_0}} \right).
\label{simpl5}
\end{equation}
The optimal choice is found by studying the variations with $z_0$  of the two functions of $\abs{\frac{z-z_0}{z_0}}$ and $ \abs{\frac{1-z_0}{z_0}}$ and the result is
\begin{equation}
z_0= \frac{z+1}{2}, \quad \norm{\toq{Z}}= \abs{\frac{z-1}{z+1}},\quad  \norm{\toq{Z}} <1 \; \Leftrightarrow \;   z \in ]0, + \infty).
\label{simpl6}
\end{equation}
The spectral radius of $\toq{\Gamma}^0 \toq{\delta L}$ being bounded by $\norm{\toq{Z}}$, the series (\ref{LS3})- (\ref{LS5}) are convergent for all $0 < z < +\infty$. However the criterion (\ref{CV2}) does not ensure convergence of the series when $z=0$, or $z<0$, or $z=+\infty$. 
\item
{\em EM-scheme:} Finally the optimal reference medium minimizing $\norm{\toq{W}}$ is found by solving the minimization problem 
\begin{equation}
\doublelow{\text{Min} \cr
z_0} \ \text{Max} \left(  \abs{\frac{z-z_0}{z+z_0}}, \abs{\frac{1-z_0}{1+z_0}} \right).
\label{simpl7}
\end{equation}
The optimal choice is found by studying the variations with $z_0$  of the two functions of $\abs{\frac{z-z_0}{z+z_0}}$ and $ \abs{\frac{1-z_0}{1+z_0}}$ and the result is
\begin{equation}
z_0= \sqrt{z}, \quad \norm{\toq{W}}= \abs{\frac{\sqrt{z}-1}{\sqrt{z}+1}}, \quad  \norm{\toq{W}} <1 \; \Leftrightarrow \; z \in ]0, + \infty).
\label{simpl8}
\end{equation}
The spectral radius of $\toq{H}^1 \toq{W}$ being bounded by $\norm{\toq{W}}$, the series (\ref{LS7})- (\ref{LS9}) are convergent for all $0 < z < +\infty$. Again the criterion (\ref{condconv2}) does not ensure convergence of the series when $z=0$, or $z<0$, or $z=+\infty$. 
\end{enumerate}
The aim of the following section is to show that, with additional information on the microstructure, it is possible to improve on the conservative estimates (\ref{simpl4}),  
 (\ref{simpl6}) and (\ref{simpl8}) for the radius of convergence of the power series.

\section{Improving the estimation of the domain of convergence of iterative methods with information on the microstructure} \label{sec3}
\subsection{What analytic properties of effective properties tell about the radius of convergence of iterative methods}
Analytic properties of the effective moduli as function of the contrast between the phases can be used to extend the domain of convergence of the power series. In the interest of simplicity we will consider here only composites following the requirements of section \ref{simpl}, where the contrast can be measured by a scalar parameter.
{Two features of analytic functions
are helpful in estimating the domain of convergence of iterative methods: }
\begin{enumerate}
\item
If $f$ is an analytic function on a domain $\mathcal D$ of the complex plane, it coincides with its Taylor series at any point of $\mathcal D$, in any disk centered at that point and lying within $\mathcal D$.
\item
If $f$ is analytic on $\mathcal D$ and $g$ is analytic on $\Omega$ containing the range of $f$ then $g(f(z))$ is analytic on $\mathcal D$.
\end{enumerate} 
\vskip 0.1cm
 Bergman \cite{BER78} proposed that in conductivity problems $\eff{\toq{L}}$ is an analytic function of $z$ in the complex plane minus the negative real axis. {A correct
argument validating this was given in Milton \cite{Milton:1981:BCP} and} a more rigorous proof was  given by Golden and Papanicolaou \cite{GOL83}. 
Therefore $\eff{\toq{L}}(z)$ can be expanded into a power series in the neighborhood of any point in the complex plane minus the negative real axis.  This result is valid for {\em any microstructure}. The radius of convergence of this series is {at least} that of the largest disk contained in the complex plane minus the real negative axis. For instance,  the expansion of $\eff{\toq{L}}(z)$ in the neighborhood of $z=1$ has a radius of convergence at least equal to $1$. 
 \vskip 0.1cm
Interestingly, the domain of analyticity of $\eff{\toq{L}}(z)$ can be extended when more information about the singularities (poles or branch-cuts) of the function $\eff{\toq{L}}(z)$ is available. This additional information often comes from an explicit expression of $\eff{\toq{L}}$ as a function of $z$. For instance, the  two Hashin-Strikman bounds on the effective conductivity of isotropic two-phase composites, which  are attained by a wide class of isotropic microstructures, have a single pole on the negative axis. Their typical form is  
\begin{equation}
\eff{\toq{L}}(z) = \eff{z} \toq{I}, \quad z^*= 1 + \alpha  \frac{z-1}{z+ \beta},
\label{HS}
\end{equation}
where $\alpha$ and $\beta$ are positive scalars depending on the volume fraction of the phases, on the dimension of space and on the type of bound (upper or lower) (see Milton \cite{MIL02}). Let us also mention that if the Hashin-Shtrikman bound (\ref{HS}) has a single isolated pole, there are other microstructures for which the whole real negative axis is singular. The checkerboard, with effective conductivity $\eff{z}=\sqrt{z}$, is an example of such a situation. In between, there are microstructures where the poles or branch-cuts are localized. 
Even in the absence of complete information about the microstructure one can sometimes restrict the interval on the real negative axis where singularities occur \cite{Bruno:1991:ECS}.
If inclusions have sharp corners then this too gives information about the singularities (Hetherington and Thorpe \cite{Hetherington:1992:CSC} and Milton \cite{MIL02}, sect.18.3).

Therefore, when the locations of the singularities  of  $\eff{\toq{L}}(z)$ are known {or partially known}
for a given microstructure, its {known} domain of analyticity may be significantly larger than the complex plane minus the negative real axis. Consequently the radius of convergence of power series can be automatically extended to all disks contained in this domain. The counterpart for Neumann series is that their radius of convergence might be significantly larger than the conservative estimates of section \ref{unconditional}. An example of such a situation will be considered in section \ref{section_theoretical}.

\subsection{Power-series expansions accounting explicitly for the microstructure}\label{sect3.2}
Apart from the analyticity of the effective moduli with respect to the contrast, it is possible to account explicitly for the microstructure by incorporating it in the operators entering the decomposition of $\toq{\Gamma}^0 \toq{\delta L}$. Instead of the decomposition (\ref{split})  one can write  
 \begin{equation}
\toq{\Gamma}^0 \toq{\delta L} =\sum_{r=1}^N  \toq{\Gamma}^1 \chi^{(r)} \toq{Z}^{(r)}, \quad \toq{Z}^{(r)} = {\toq{L}^0}^{-1} \toq{\delta L}^{(r)},
\label{iter3}
\end{equation}
where $\chi^{(r)}$ is the characteristic function of phase $r$. The operator $\toq{\Gamma}^0 \toq{\delta L} $ entering the Liouville-Neumann series (\ref{LS3})-(\ref{LS5})  is therefore the composition of two operators: 
\begin{enumerate}
\item
The operators $\toq{Z}^{(r)}$ depend only the {\em material properties} of the phases (and on the reference medium). They serve as measures of the contrast of the composite.
\item
The operators $\toq{\Gamma}^1 \chi^{(r)}$ depend on the {\em microstructure} of the composite through $\chi^{(r)}$ (but not exclusively since $\toq{\Gamma}^1$ may depend on the reference medium). These operators express the way in which the operator $\toq{\Gamma}^1$ ``sees'' the microstructure.
\end{enumerate}
The decomposition (\ref{iter3}) highlights  the contribution of the material properties on one hand and of the microstructure on the other hand. In a quite similar way the operator entering the Eyre-Milton scheme can be written as 
\begin{equation}
\toq{H}^1 \toq{W} =\sum_{r=1}^N  \toq{H}^1 \chi^{(r)} \toq{W}^{(r)}, \quad \toq{W}^{(r)} = (\toq{L}^{(r)}+\toq{L}^0)^{-1} \toq{\delta L}^{(r)}.
\label{iter3bis}
\end{equation}
The general case being rather technical, attention is limited here to two-phase composite whose contrast can be measured by a single scalar $z$, as described in section \ref{simpl}. The expression of $\toq{Z}$ and $\toq{W}$ are given in (\ref{simpl2})
and the decompositions (\ref{iter3}) and (\ref{iter3bis}) read as 
\begin{equation}
\left.
\begin{array}{c}
\displaystyle \toq{\Gamma}^0 \toq{\delta L}=\toq{\Gamma}^1 \toq{Z} =\frac{z-z_0}{z_0} \Gc+ \frac{1-z_0}{z_0}\Guc,\\[2ex] 
\displaystyle \toq{H}^1 \toq{W}=
\frac{z-z_0}{z+z_0} \Hc+ \frac{1-z_0}{1+z_0}\Huc,
\end{array}
\right\}
\label{iter4}
\end{equation}
where  
\begin{equation}
\left.
\begin{array}{c}
\Gcr= \toq{\Gamma}^1 \chi^{(r)} \tod{I}, \quad \Gc+\Guc=\G,\quad \G= \toq{\Gamma}^1\tod{I}, \\ [2ex]
\Hc= \toq{H}^1 \chi^{(1)} \tod{I},   \quad \Hc+\Huc=\Hh,\quad \Hh= \toq{H}^1\tod{I}. 
\end{array}
\right\}
 \label{iter4bis}
 \end{equation}
It follows from (\ref{iter4}) that the key operators in the iterations for the three cases of interest identified in section \ref{simpl} are: 
\vskip 0.1cm
\noindent
{\em Reference medium= Matrix (B-scheme):}
\begin{equation}
 \left.
\begin{array}{l}
z_0=1,\quad \toq{\Gamma}^0 \toq{\delta L}=t \Gc, \quad t=z-1\\[2ex]  
\displaystyle {\toq{L}^0}^{-1} \eff{\toq{L}}= \toq{I} + \sum_{j=0}^{\infty} (-1)^j \moy{\chi^{(1)}( \Gc)^{j}}t^{j+1}.
\end{array}
\right\}
\label{recap1}
\end{equation}
{\em Reference medium= Arithmetic mean (MS-scheme):}
\begin{equation}
\left.
\begin{array}{l}
\displaystyle z_0=\frac{z+1}{2}, \quad
\toq{\Gamma}^0 \toq{\delta L} = t \left(\Gc - \Guc\right), \quad t= \frac{z-1}{z+1}, \\[2ex]
{\displaystyle {\toq{L}^0}^{-1} \eff{\toq{L}}= \toq{I} + \sum_{j=0}^{\infty} (-1)^j  \moy{(\chi^{(1)}-\chi^{(2)}) (\Gc-\Guc)^{j}} t^{j+1}.} 
\end{array}
\right\}
\label{recap2}
\end{equation}
{\em Reference medium= Geometric mean (EM-scheme):}
\begin{equation}
\left.
\begin{array}{l}
\displaystyle z_0=\sqrt{z},
\toq{H}^1 \toq{W} = t \left(\Hc - \Huc\right), \quad t= \frac{\sqrt{z}-1}{\sqrt{z}+1}, \\[2ex]
{\displaystyle {\toq{L}^0}^{-1} \eff{\toq{L}}= \toq{I} + 2 \sum_{j=0}^{\infty} (-1)^j  \moy{(\chi^{(1)}-\chi^{(2)})  (\Hc- \Huc)^{j}} t^{j+1}.}
\end{array}
\right\}
\label{recap3}
\end{equation}
The expressions of the effective moduli in  (\ref{recap2}) and (\ref{recap3}) evidence the symmetric role played by the two phases (they also play a symmetric role in the   reference medium), whereas the first expression is asymmetric (because the reference medium is one of the two phases). The characteristic function of one phase can be eliminated from the last two relations in (\ref{recap2}) by noting that 
\begin{equation} 
\chi^{(1)}- \chi^{(2)} =2\chi^{(1)} -1\stackrel{def}{=}\chi', \quad  \Gc- \Guc= \Gcp,\quad  \Hc- \Huc= \Hcp.
\label{chiprime}
\end{equation}
{It is remarkable that, in the series expansion (\ref{recap1}),  (\ref{recap2}) and (\ref{recap3}), the variable $t$ measuring the contrast between the phases and the microstructure are decoupled.} In other words, one could compute  the micro-structural coefficients $\moy{\chi^{(1)}( \Gc)^{j}}$,  $\moy{\chi' (\Gcp)^{j}}$, $\moy{\chi' (\Hcp)^{j}}$ once for all. And for any arbitrary contrast $t$ the effective property would be obtained by summing the series (\ref{recap1}),  (\ref{recap2}) or (\ref{recap3}), without having to solve the Lippmann-Schwinger equation. {A similar observation has previously been made by Hoang and Bonnet \cite{HOA14} for incompressible phases and recently by To {\em et al}\cite{TO17}  for elastic phases with the same Poisson ratio. It is also closely related to the role of correlation functions of increasing order into the series expansions of the effective moduli (Milton \cite{MIL02} chap. 15).}

\section{A study case: the {Obnosov} formula for the conductivity of a square array of square inclusions}
\label{section_theoretical}
\begin{figure}[!h]
  \begin{center}
  \parbox{5cm}{
    \includegraphics[width=5cm,height=5.cm]{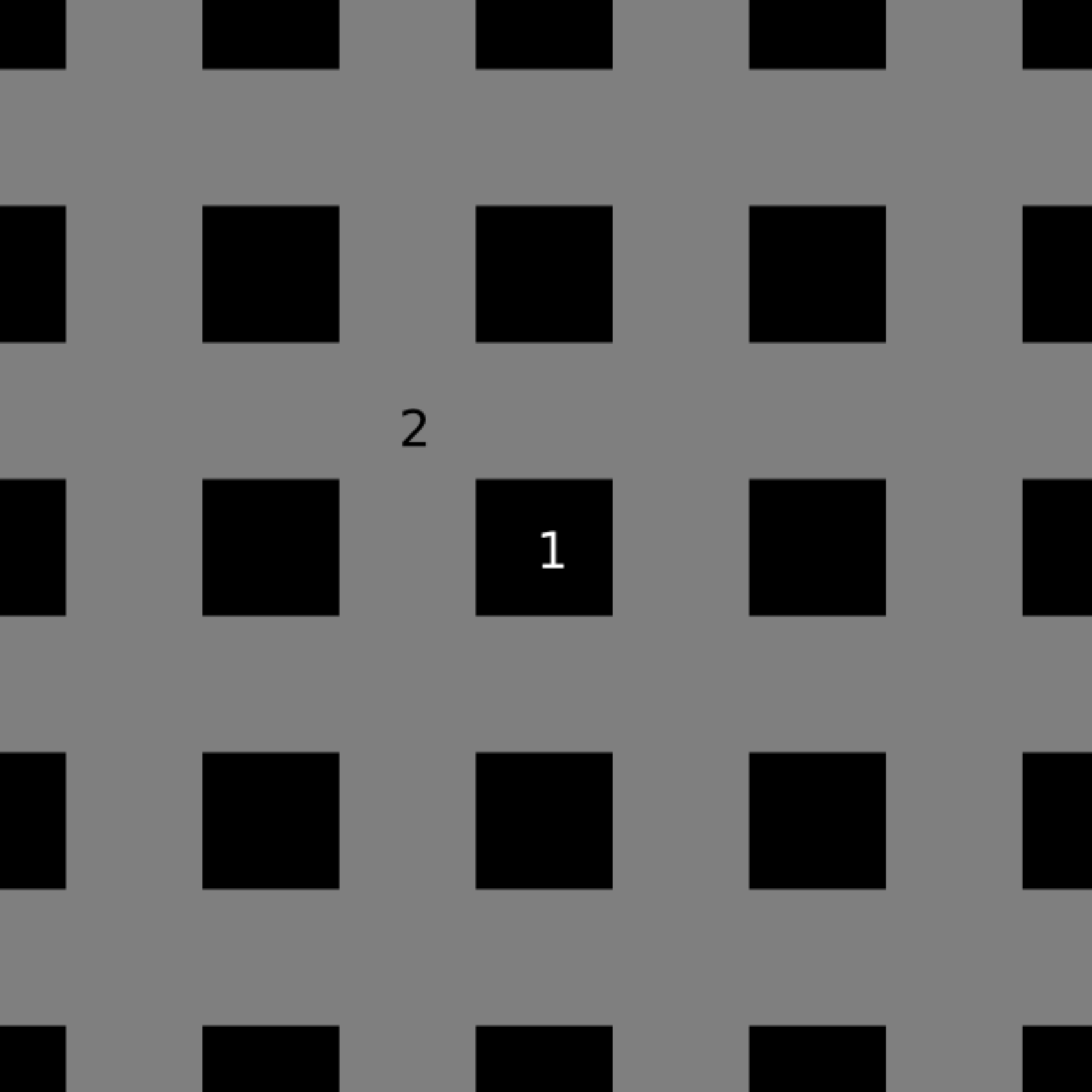}
  }
  \end{center}
  \caption{Obnosov microstructure: square array of square inclusions with volume fraction $0.25$.}
  \label{Image}
\end{figure}

Obnosov \cite{OBN99} derived a closed form expression for the {effective} conductivity of a square array of square {isotropic} inclusions with volume fraction $0.25$ (cf Figure \ref{Image}). 
This microstructure has square symmetry and its effective conductivity tensor is isotropic
\begin{equation}
\eff{\toq{L}}(z) = \eff{z} \toq{I}, \quad  \eff{z} = \sqrt{   \frac{ 1 + 3 {z} }{ 3 + {z}}   },
  \label{obnosov1}
\end{equation}
where, after proper rescaling the conductivity of phase 1 (inclusion) is $z$ and that of phase 2 (matrix) is $1$. 

\subsection{Convergence of the exact power series}
Exact power series expansions for $\eff{z}${, or alternatively for  $\eff{z}/z_0$} in the form
\begin{equation}
  \eff{z} = \sum_{k=0}^\infty  {d_k} t^k,\quad {\frac{\eff{z}}{z_0}=\sum_{k=0}^\infty  {b_k} t^k,}
  \label{series_scalar}
\end{equation}
can be derived from (\ref{obnosov1}) in the three cases of interest described in section \ref{simpl}, where the contrast variable $t$ {depends} on the choice of the reference medium.

\begin{figure}[htp!]
 \parbox{5.2cm}{
\includegraphics[width=5.5cm]{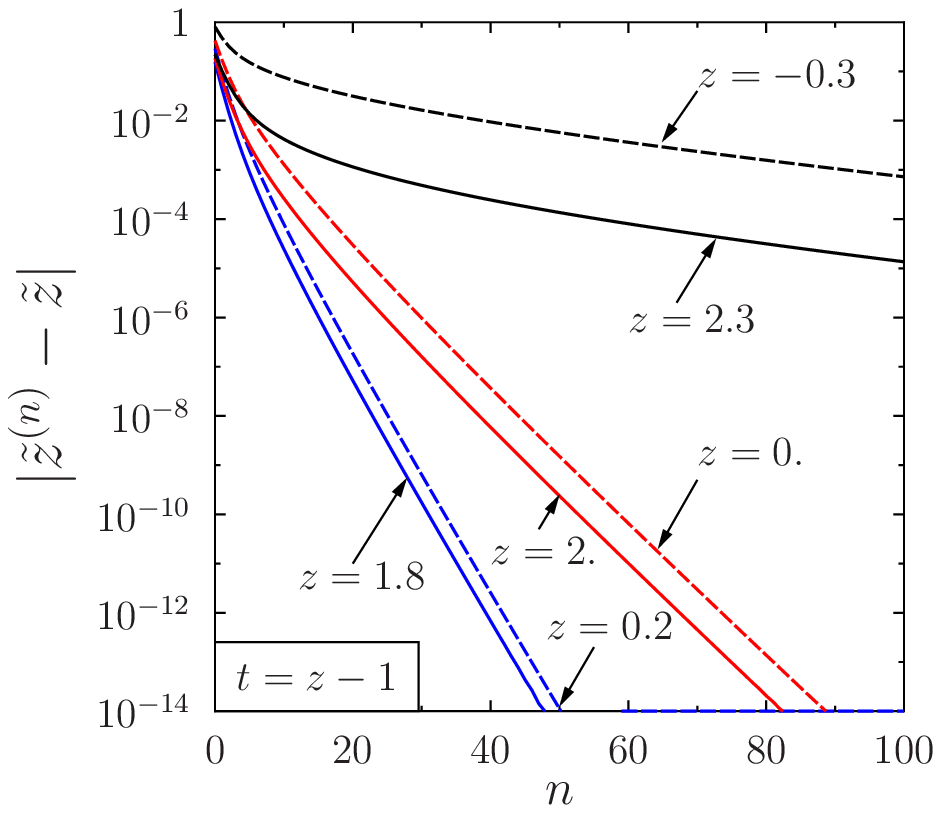}
\begin{center}
\vskip -0.5cm
$\qquad $  (a) 
\end{center}    }
     \parbox{5.2cm}{
      \includegraphics[width=5.5cm]{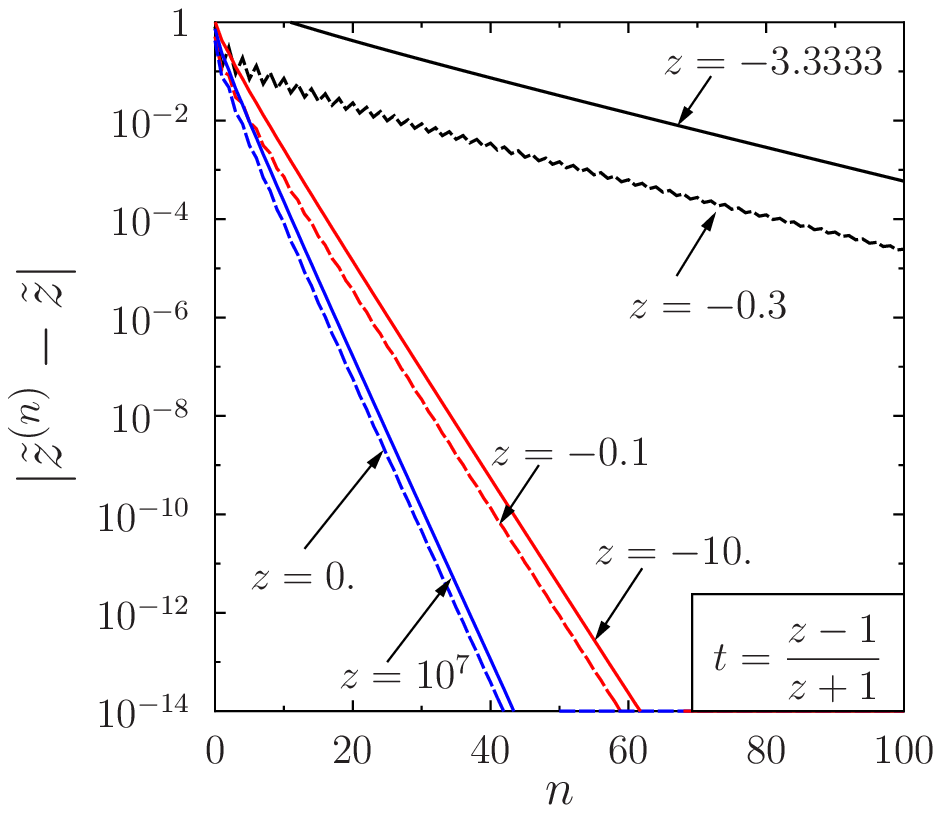}
 \begin{center}
\vskip -0.5cm
 $\qquad $  (b)
\end{center}      }
    \parbox{5.2cm}{
          \includegraphics[width=5.5cm]{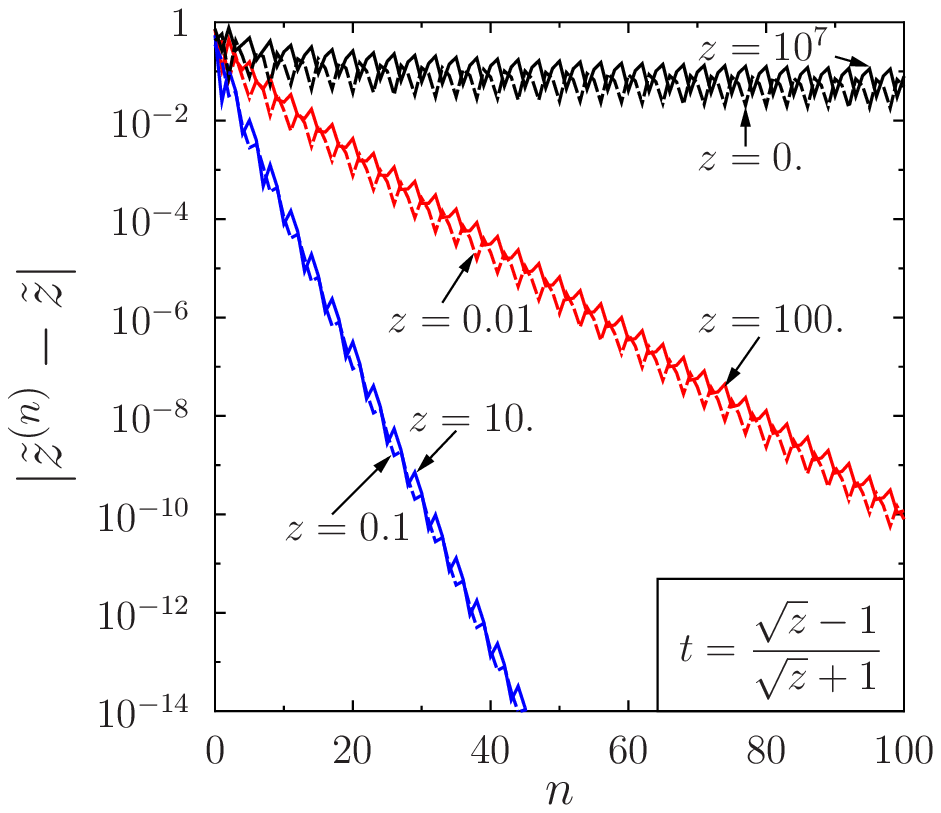}
 \begin{center}
\vskip -0.5cm
$\qquad $ 
(c)
\end{center}      }
    \caption{Convergence of the theoretical power series (\ref{series_scalar}). Error (\ref{err}) as a function of the order of truncation.
{Expansion with respect to (a)
 $t=z-1$ , (b) $ \displaystyle t = \frac{z-1}{z+1}$,
 and (c) $\displaystyle t = \frac{\sqrt{z}-1}{\sqrt{z}+1}$ .
}
Solid lines: $t>0$, dashed lines $t<0$. }
  \label{fig_zm1}
\end{figure}

Typically, $t=z-1$ or {$\displaystyle t= \frac{z-1}{z+1}$} or 
{$\displaystyle  t=\frac{\sqrt{z}-1}{\sqrt{z}+1}$.} These power series 
are given in appendix \ref{dvpt} for the three cases of interest. 
\vskip 0.1cm
The conservative estimates of section \ref{unconditional} ensure convergence of these power series for $z \in \mathbb{C} - \mathbb{R}^-$, {\em i.e.} when $z$ satisfies the conditions (\ref{simpl4}) or (\ref{simpl6}) (which coincides with (\ref{simpl8})).  These conservative estimates can be improved by means of the analytic properties of the effective conductivity (\ref{obnosov1}). 
\vskip 0.1cm
Consider first the expansion with respect to $z-1$ (corresponding to the matrix as reference medium).  $\eff{z}(z)$  has a branch-cut $[-3,-\frac{1}{3}]$ on the negative real axis but is analytic in the whole complex plane outside this branch-cut. Therefore the power series of $\eff{z}$ in the neighborhood of  $z=1$ will converge in all disks centered at $z=1$ and not intersecting the branch-cut, {\em i.e.} for $\abs{z-1} < {4}/{3}$. Considering only real values for  $z$, the power series converges when $- {1}/{3} < z < {7}/{3}$ and in particular when $z=0$. This extended domain of convergence is confirmed in Figure \ref{fig_zm1}(a) where the deviation from the exact result (\ref{obnosov1}) as a function of the number $n$ of terms in the series is measured by 
\begin{equation}
 {\epsilon}_n = \abs{ \eff{z}^{(n)} - \sqrt{\frac{1+3z}{3+z} }}, \quad \eff{z}^{(n)} = \sum_{k=0}^n d_k t^k.
\label{err}
\end{equation}
Note in particular that the series converges (although slowly) even when $z=-0.3$, in accordance with the theoretical predictions.
\vskip 0.1cm
The situation is similar for the power series in {$t = \frac{z-1}{z+1}$}. The effective conductivity (\ref{obnosov1}) reads as 
{
  \begin{equation}
\eff{z} = \sqrt{\frac{1+t/2}{1-t/2}},
\label{obnosov2}
  \end{equation}
  }
which has two branch-cuts $(-\infty,-2]\cup[2,+\infty)$. The disk of convergence of its power series expansion in the neighborhood of {$t=0$} is therefore {$\abs{t}<2$}. 
\vskip 0.1cm
Limiting attention to $z$ on the real axis, the power series in {$t= \frac{z-1}{z+1}$} converges when 
$z \in (-\infty,-3[\cup]-1/3,+\infty)$. Therefore the range of $z$ for which this series converges is larger than in the previous expansion with respect to $z-1$. This theoretical prediction is  confirmed by the numerical experiments reported in figure \ref{fig_zm1}(b) where the error between the partial series and the exact effective conductivity  is still measured by (\ref{err}), but with $t= \frac{z-1}{z+1}$. Note in particular that the series converges for $z=-10$.
\vskip 0.1cm
The expansion with respect to $t= {(\sqrt{z}-1)}/{(\sqrt{z}+1)}$ is slightly different since, in order for $\sqrt{z}$ to be properly defined, $z$ has to be in the cut complex plane (outside the negative real axis).


The convergence of the corresponding series is quite fast at moderate contrast, but is not very fast at high contrast ($z=0$ and $z=10^7$). {This is expected as the square root in $\sqrt{z}$
introduces a branch cut of singularities extending from zero to infinity.}
\subsection{Convergence of numerical series}  
As already mentioned in section {\ref{sect3.2} the power series expansion of the effective tensor $\eff{\toq{L}}$ can be obtained by successive applications of the operators $\toq{\Gamma}^1$ or $\toq{H}^1$. 
Taking into account the specific form of the effective tensor $\eff{\toq{L}}= \eff{z} {\toq{I}}$ and  $\toq{L}= {z}_0 {\toq{I}}$ into (\ref{recap1}),  (\ref{recap2}) and (\ref{recap3}) the coefficients $d_k$ of the power series (\ref{series_scalar}) can be related to the characteristic function of the phases and to  the operators   $\toq{\Gamma}^1$ or $\toq{H}^1$ through the following relations: 
\vskip 0.1cm
\noindent
{\em Reference medium= Matrix (B-scheme):} $t=z-1$  
\begin{equation}
d_0=1, \; d_k= {b_k}= (-1)^{k-1} \moy{\chi^{(1)}( \Gc)^{k-1}}.\tou{e}_1 \quad k \geq 1.
\label{dkB}
\end{equation}
{\em Reference medium= Arithmetic mean (MS-scheme):} $t= \displaystyle  \frac{z-1}{z+1}$
\begin{equation}
d_0=1, \; d_k-d_{k-1} = {b_k, \; b_k=}(-1)^{k-1}  \moy{\chi' (\Gcp)^{k-1}}.\tou{e}_1,\quad k \geq 1.
\label{dkMS}
\end{equation}
{\em Reference medium= Geometric mean (EM-scheme):} $t= \displaystyle \frac{\sqrt{z}-1}{\sqrt{z}+1}$
\begin{equation}
{d_0=1,\; d_1=2 (\moy{\chi'} +1),\;  d_k-d_{k-1} = {b_k + b_{k-1},\; b_k=2 (-1)^{k-1}   \moy{\chi'  (\Hcp)^{k-1}}.\tou{e}_1,}\quad k \geq 2,}
\label{dkEM}
\end{equation}
where $\chi'$, $\Gcp$ and $\Hcp$ are defined in (\ref{chiprime}). The $\Gc$, $\Gcp$ and $\Hcp$  and the operator $\toq{\Gamma}^0$ from which they are formed, are more easily handled in Fourier space where they are explicitly known. {In conductivity and elasticity the operator }
$\toq{\Gamma}^0$ can be expressed in Fourier space as 
\begin{equation}
\hat{\toq{\Gamma}}^0(\tou{\xi})= \tou{\xi} \otimes \left(  \tou{\xi}.\toq{L}^0. \tou{\xi} \right)^{-1} \otimes  \tou{\xi}.
\label{Gamma0}
\end{equation}
In the above expression of the Green's operator,  $\tou{\xi}$ is the {spatial frequency} in Fourier space. 
Therefore, in practice, the coefficients of the series expansion (\ref{series_scalar})  can be computed by successive applications of the same operator $\Gc$, $\Gcp$ or $\Hcp$ depending on the reference medium.
Each application of this operator involves a multiplication by a characteristic function followed by a Fourier transform. The latter operation is performed on a discretized image of the microstructure {by Fast Fourier Transform}. Therefore computing $d_k$ will require, in all three cases, $k-1$ {applications} of the FFT to a non{smooth} field (the non{smoothness} 
stemming from the characteristic function). 
The ``theoretical power series'' (TS) (\ref{series_scalar}) (where the exact coefficients {$b_k$ and $d_k$} are given {for each scheme} in appendix \ref{dvpt})  and the  ``numerical power series'' (NS) (\ref{dkB}, \ref{dkMS}, \ref{dkEM}) are compared in figure \ref{comparaison1} for the three schemes of section \ref{simpl}. 
\begin{figure}[htp!]
    \parbox{5.2cm}{
      \includegraphics[width=5.5cm]{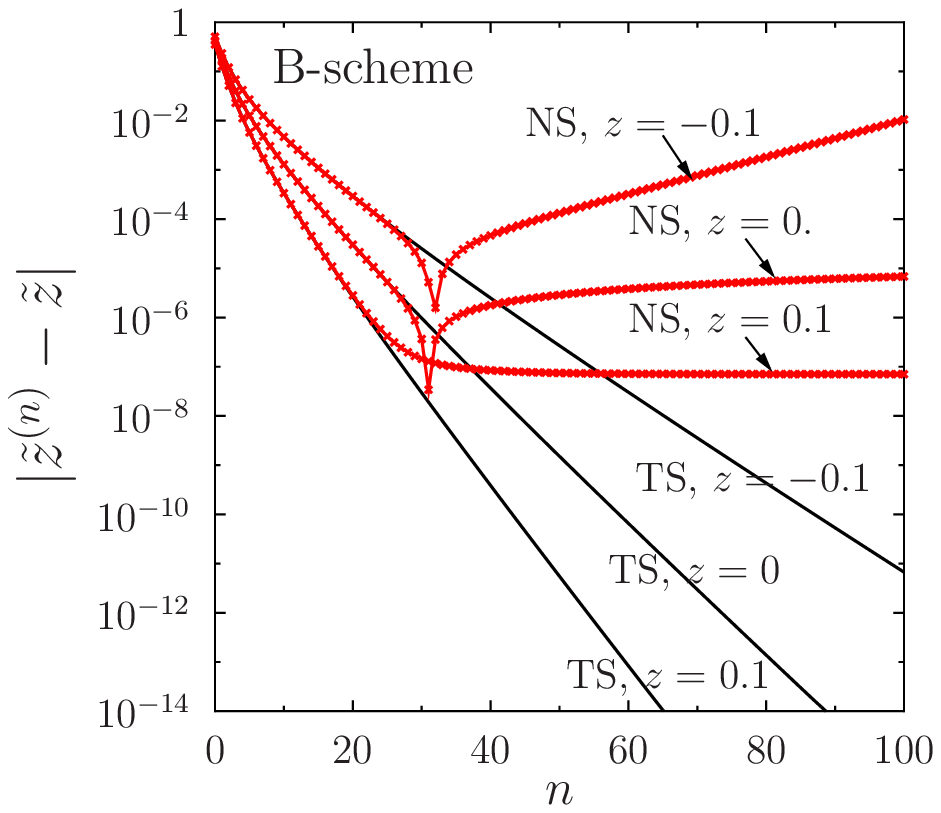}       
      \begin{center} \vskip -0.5cm   $\qquad$ (a)  \end{center}
    }
    \parbox{5.2cm}{
      \includegraphics[width=5.5cm]{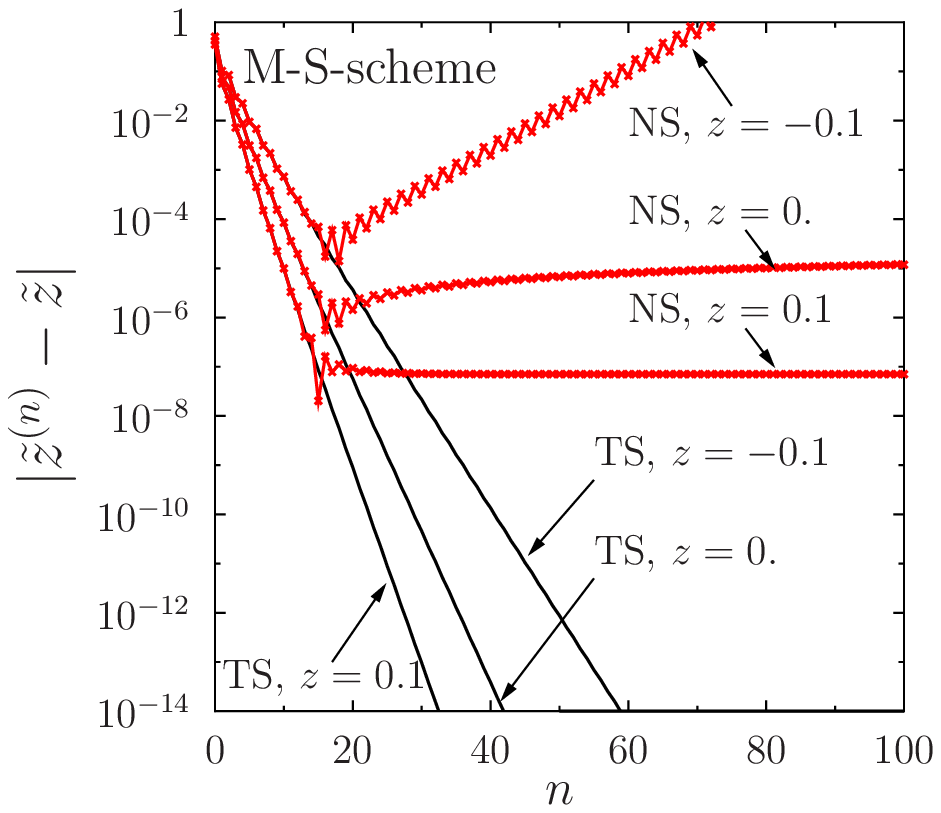}      
     \begin{center} \vskip -0.5cm  $\qquad$ (b) \end{center}
    }
        \parbox{5.2cm}{
     \includegraphics[width=5.5cm]{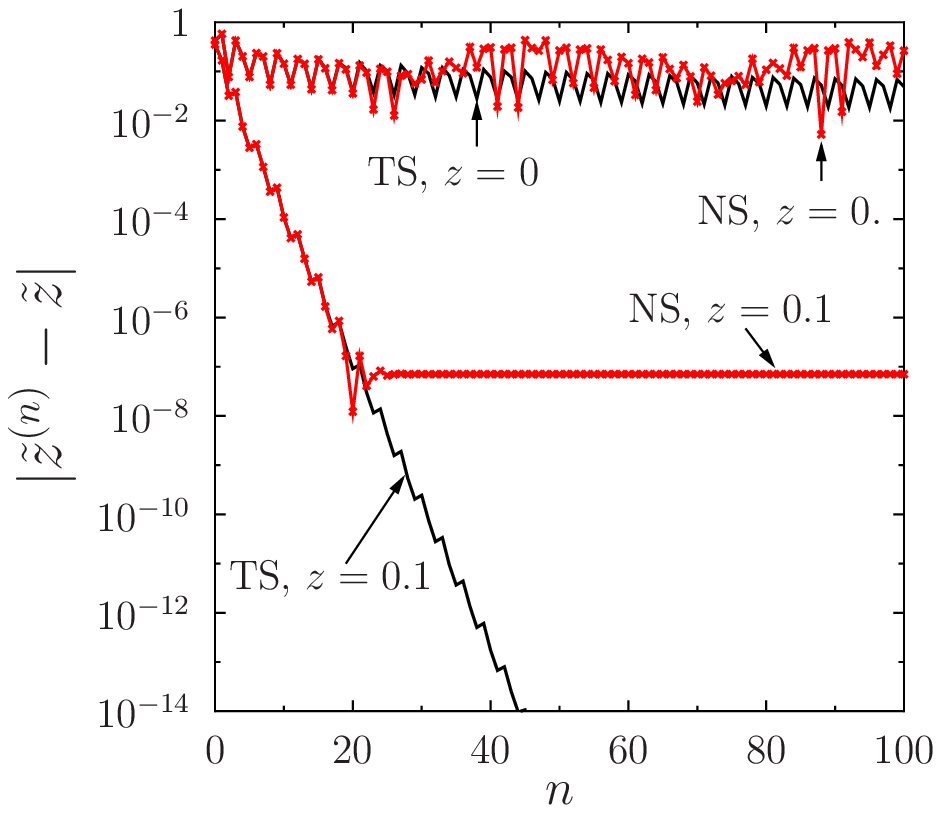}      
        \begin{center} \vskip -0.5cm $\qquad$ (c) \end{center}
    }
    \caption{Obnosov problem. Comparison between the convergence of the theoretical series (TS) and the numerical  series (NS) for different contrasts. (a) Reference medium: matrix (B-scheme). (b) Reference medium: arithmetic mean (MS-scheme). (c) Reference medium: geometric mean (EM-scheme). Discretization:  $512 \times 512$ pixels.
 }
 \label{comparaison1}
\end{figure}
The following observations can be made:
\begin{enumerate}
\item
For all three contrast variables, the sums of the first 10 to 20 terms of the theoretical and numerical series coincide. 
\item
However, after a certain number of iterations, the deviation of the computed effective conductivity from the exact result decreases only slowly (for positive contrast) or saturates, or even increases (for negative contrast) when more and more terms are added.
\item
Figure \ref{comparaison2} shows that the point of deviation between the two series depends on the discretization of the image. Three discretization for the image, $128^2$, $512^2$ and $2048^2$ pixels respectively, have been considered. The finer the discretization, the later the deviation and the better the prediction of effective conductivity. In the limit of an infinitely fine discretization, it is expected that both series will coincide at any order. 
\end{enumerate}
\begin{figure}[!]
    \parbox{5.2cm}{
            \includegraphics[width=5.5cm]{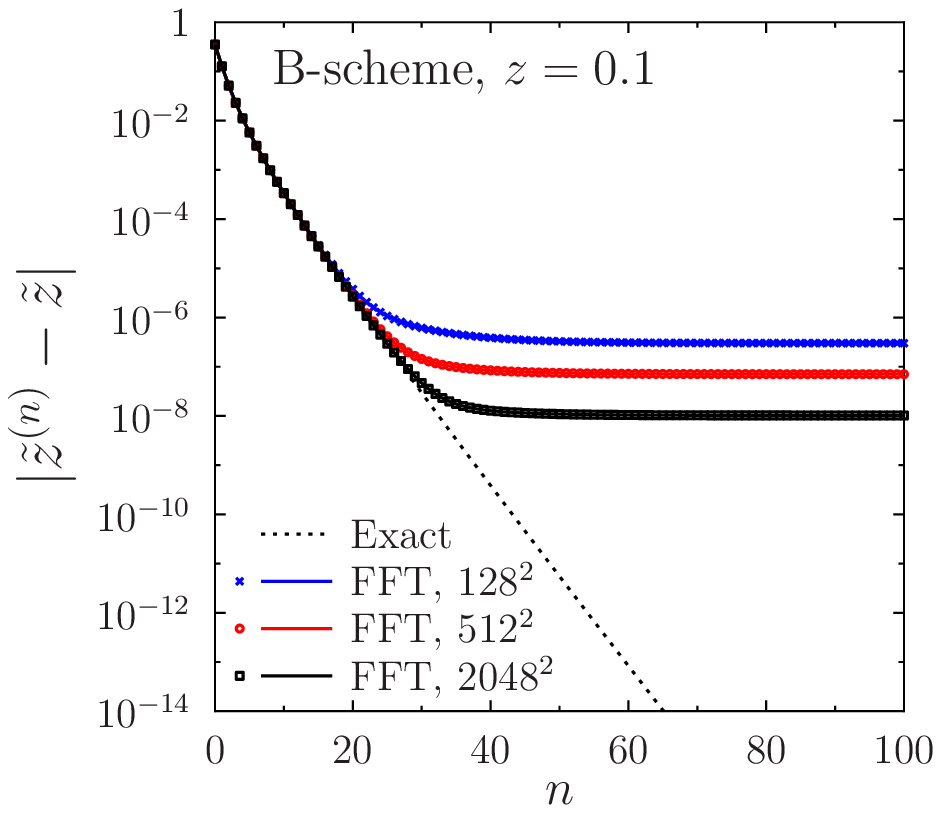}
          \begin{center} \vskip -0.5cm
$\qquad$      (a) \end{center}
    }
    \parbox{5.2cm}{
         \includegraphics[width=5.5cm]{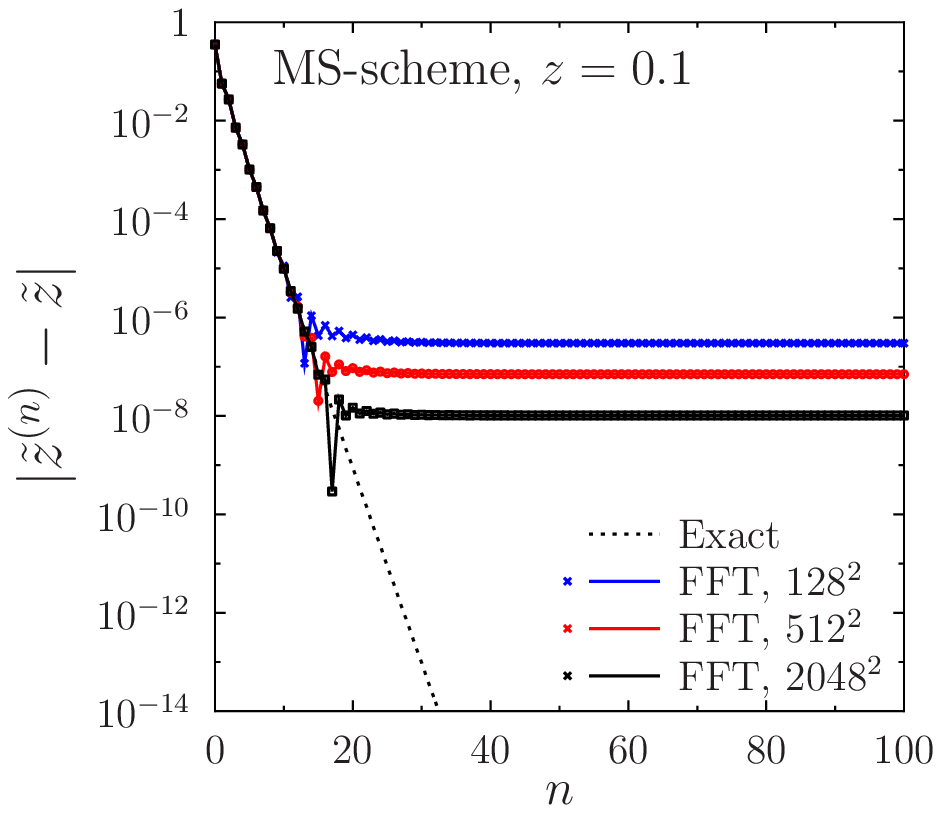}     
                   \begin{center}  \vskip -0.5cm  $\qquad$ (b) \end{center}
    }
    \parbox{5.2cm}{
        \includegraphics[width=5.5cm]{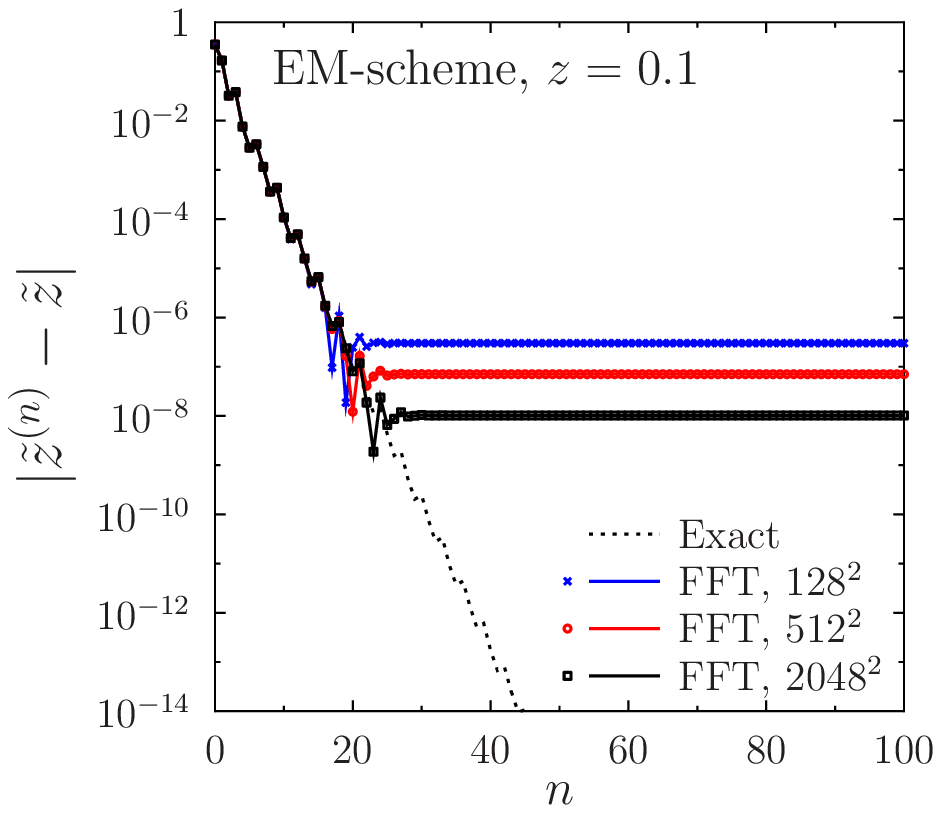}      
            \begin{center} \vskip -0.5cm   $\qquad $ (c) \end{center}
    }
    \vskip 0.1cm
        \parbox{5.2cm}{
            \includegraphics[width=5.5cm]{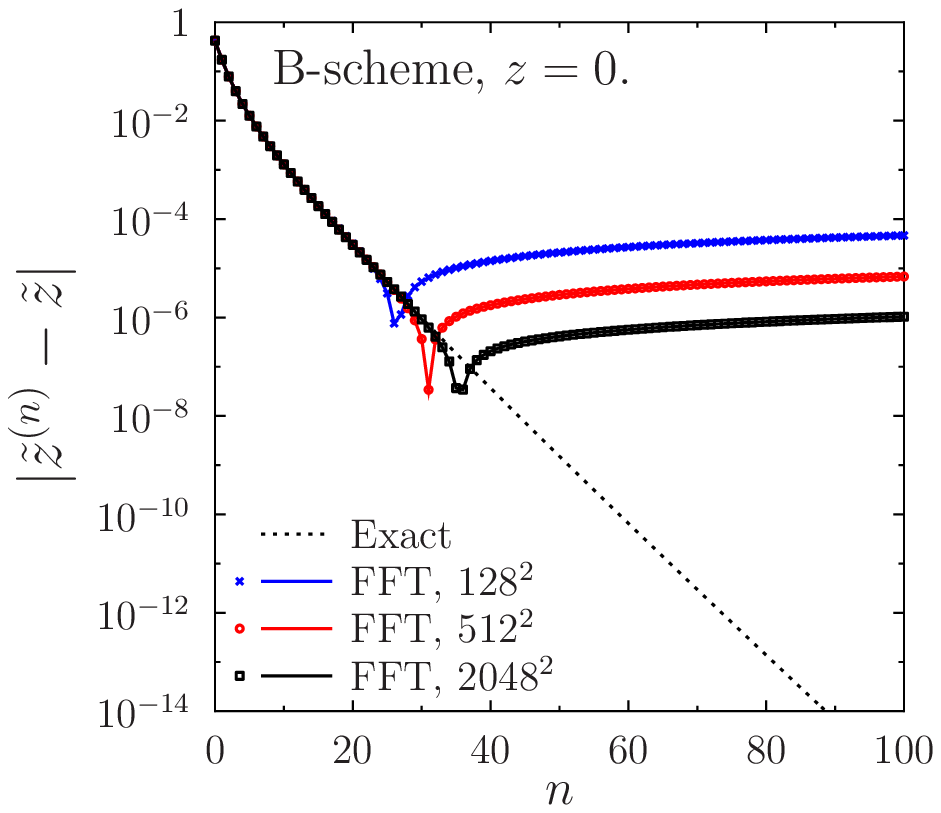}
          \begin{center} \vskip -0.5cm
    $\qquad$  (d) \end{center}
    }
    \parbox{5.2cm}{
         \includegraphics[width=5.5cm]{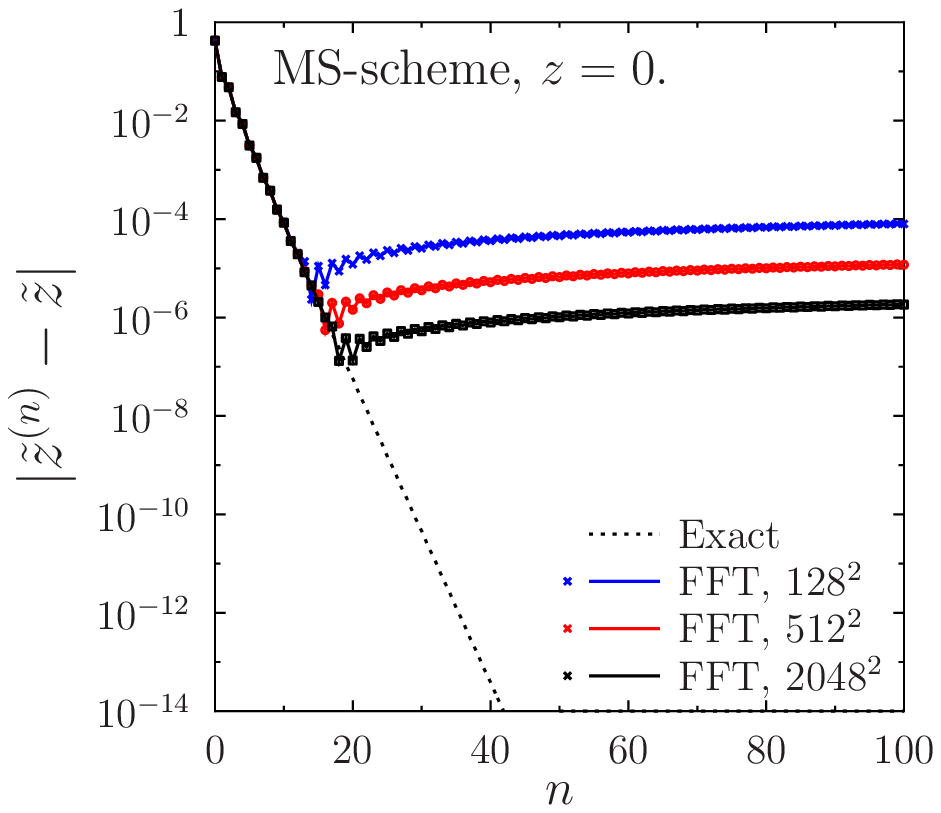}     
                   \begin{center}  \vskip -0.5cm  $\qquad$ (e) \end{center}
    }
    \parbox{5.2cm}{
        \includegraphics[width=5.5cm]{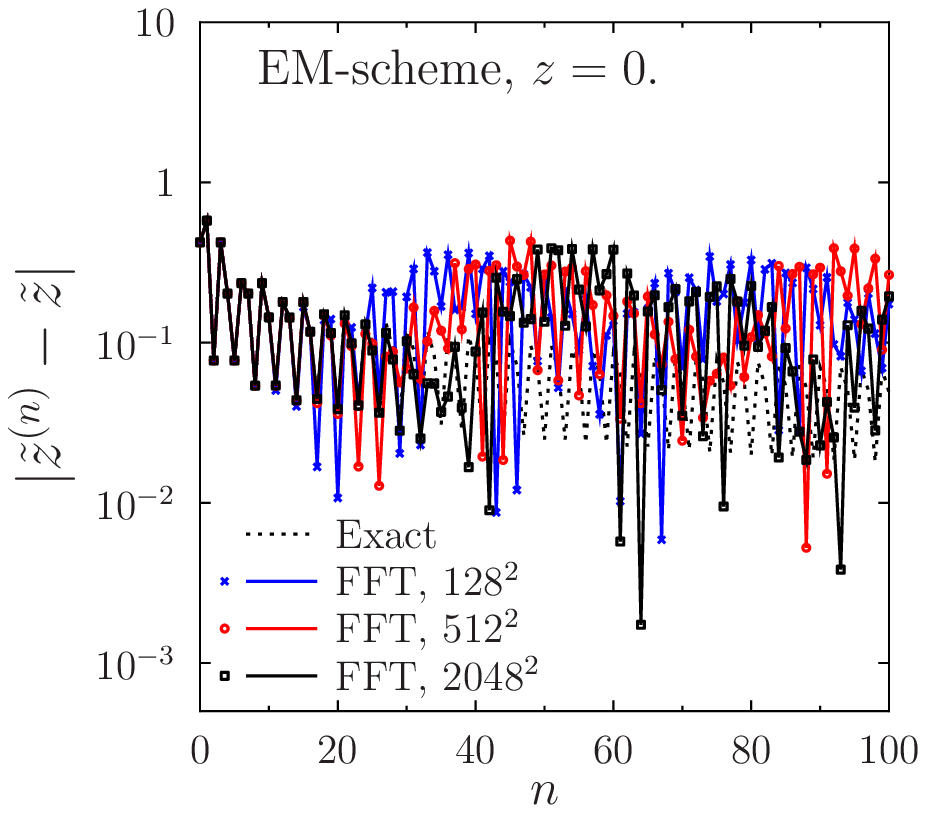}      
            \begin{center} \vskip -0.5cm   $\qquad$ (f) \end{center}
    }
    \caption{Influence of the discretization on the convergence of the numerical scheme. (a)(b)(c): contrast $z=0.1$, (d)(e)(f): $z=0$. (a)(d) Reference medium: matrix (B-scheme). (b)(e) Reference medium: arithmetic mean (MS-scheme). (c)(f) Reference medium: geometric mean (EM-scheme). }
  \label{comparaison2}
\end{figure}
\subsection{Discussion}\label{discussion}
{
Even with the discretized equations, the truncated series expansions are rational functions of the {contrast} and therefore in the domain of convergence they converge to an analytic
function (\cite{Rudin:1987:RCA}, theorem 10.28). The analytic continuation of this function will have singularities at the radius of convergence. Such singularities might include, for example, a pole with
very small residue that is not necessarily located on the negative real axis. These spurious singularities in the discrete approximation dictate the convergence properties of
the discrete approximation.
}

The only difference between the theoretical and the numerical series $\eff{z}^{(n)}$ resides in the coefficients $d_k$.
The numerical coefficients, which will be denoted $d^{NS}_k$ to distinguish them from the exact ones $d_k$, are only an  approximation of the exact coefficients, the accuracy of the approximation depending on the refinement of the discretization and on the choice of the $\toq{\Gamma}^0$ operator (continuous or discrete). 
\vskip 0.1cm
The role of the coefficients $d_k$ in the discrepancy between the theoretical and numerical predictions of the effective conductivity is confirmed by figure {\ref{comparaison3}} and table \ref{table_dk} where the first $d_k$ coefficients in the theoretical and numerical series
are compared.    It is seen by comparing  figure \ref{comparaison1} and figure {\ref{comparaison3}} that the deviation in the effective properties observed in figure \ref{comparaison1}  coincides with a similar deviation in the coefficients. Table \ref{table_dk} shows the variation of the individual coefficients with the order in the series (or with the number of iterations). 
It clearly appears that beyond a certain order in the expansion (or beyond a certain number of iterations), hereafter denoted by $\mathcal{K}$,
there is a strong discrepancy between the theoretical and numerical $d_k$'s.
The first terms (up to  {$\mathcal{K}$}) in the numerical series can be considered as exact, whereas the terms of higher order differ from the exact values. Therefore the series expansion for $\eff{z}$ may be written as 
\begin{equation}
\eff{z}^{NS}_n \simeq \sum_{k=0}^\mathcal{K} {d}_k t^k +
  \sum_{k=\mathcal{K}+1}^n d^{NS}_k t^k \; .
  \label{Dk}
\end{equation}
The coefficients $d^{NS}_k$ terms ($k>\mathcal{K}$)  decrease slowly (at best, or even remain constant) when $k$ increases but remain much larger than the exact $d_k$ terms. This discrepancy is indeed not surprising, since if all coefficients were accurately computed at any order, the sum of the series would be the {\em exact} effective property. However, the use of any computational method means that the solution of the continuous problem, which belongs in general to a vector space of infinite dimension, is approximated in a space of finite dimension. Therefore an intrinsic error due to the spatial discretization exists, resulting in a gap between the exact solution and its finite dimensional approximation, regardless of the computational method.   And since the effective property is computed here by means of its variational characterization (through the associated energy), this gap is always non zero regardless of the computational method (except when the solution belongs to the finite dimensional space of approximation, which is rarely the case). Therefore the sum of the approximate series differs from the exact solution and the coefficients of the numerical series have to deviate, sooner or later, from their theoretical value. Clearly enough, when the spatial discretization of the microstructure is refined, the approximate solution is a better approximation of the exact one, the effective energy is lower and the gap reduces. Therefore the deviation from the theoretical values is delayed by refining the discretization (but occurs eventually).   This is confirmed in figure {\ref{comparaison3}}.
\vskip 0.1cm
Our interpretation of this deviation from exactness goes as follows. In the course of the iterations, the successive approximations (partial sums of the series) converge towards the ``best'' approximation of the solution in the finite dimensional space associated with a given discretization. When the iterates are close to this best approximation, the effective energy cannot be further lowered and the effective properties reach a plateau as can be seen in figure \ref{comparaison2} (a)-(c) for instance. Iterating further neither improves the plateau, which is an energetic barrier, nor the local solution in this finite dimensional space. In conclusion, all iterations beyond $\mathcal K$ introduced in (\ref{Dk}) do not improve convergence and can even deteriorate the local field (as observed sometimes). It is therefore proposed in section \ref{newconv} to stop iterating at
{$\mathcal{K}$}, provided {$\mathcal{K}$} can be properly defined. 
The effect of the discretization on the loss of accuracy along the iterations is not restricted to the specific microstructure investigated here, it is indeed a generic observation. To illustrate this point, the MS scheme has been applied to the 2D microstructure of figure {\ref{dk_4c}(a)} and the $d_k$ terms of the series {(\ref{recap2})} have been computed numerically for different resolutions, from $64 \times 64$ pixels to $4096 \times 4096$ pixels. As for the Obnosov microstructure, the $d_k$ terms begin to deviate from each other after a certain number of iterations which depends on the discretization. The point of deviation between a coarse and a fine discretization is probably the point where the iterations for the coarser discretization should be stopped, since all iterations beyond that point are affected by a significant discretization error.
\vskip 0.1cm
The deviation from exactness of the coefficients of the series expansions is essential to understand why iterative methods fail to converge when the contrast between the constituents is very large or negative, although theoretically they should converge for certain microstructures (the Obnosov microstructure is one of them) even when one of the phases has a nonpositive modulus in a certain range. In such situations of extreme contrast, the absolute value of the parameter $t$ measuring the contrast  is close to 1 or even larger than 1.  Therefore the convergence of the series $d_k t^k$ relies crucially on the accuracy of the $d_k$'s. But since these coefficients must deviate from their exact value beyond a given order $\mathcal{K}$ it is not surprising that  the series does not converge towards its exact sum, or does not converge at all. Conversely, when the contrast parameter $t$ is strictly less than 1 (in absolute value), the deviation from exactness of the coefficients at high order does not really influence the convergence of the series $d_k t^k$ which is governed by $t^k$. This is the situation for moderate contrast between the constituents. 
\vskip 0.1cm
The threshold $\mathcal{K}$ depends on the computational method used to perform the iterations, or equivalently to compute the coefficients of the series expansions with the help of the relations of section \ref{sect3.2}. The Fast Fourier Transform, in discrete form, is  used in the present study. Its role in the deviation of the coefficients remains to quantified. {However it  is clear, qualitatively, that the deviation} will occur in general {with} {\em any} other computational method operating within a finite dimensional space. }

\begin{figure}[!]
   \parbox{5.2cm}{
      \includegraphics[width=5.5cm]{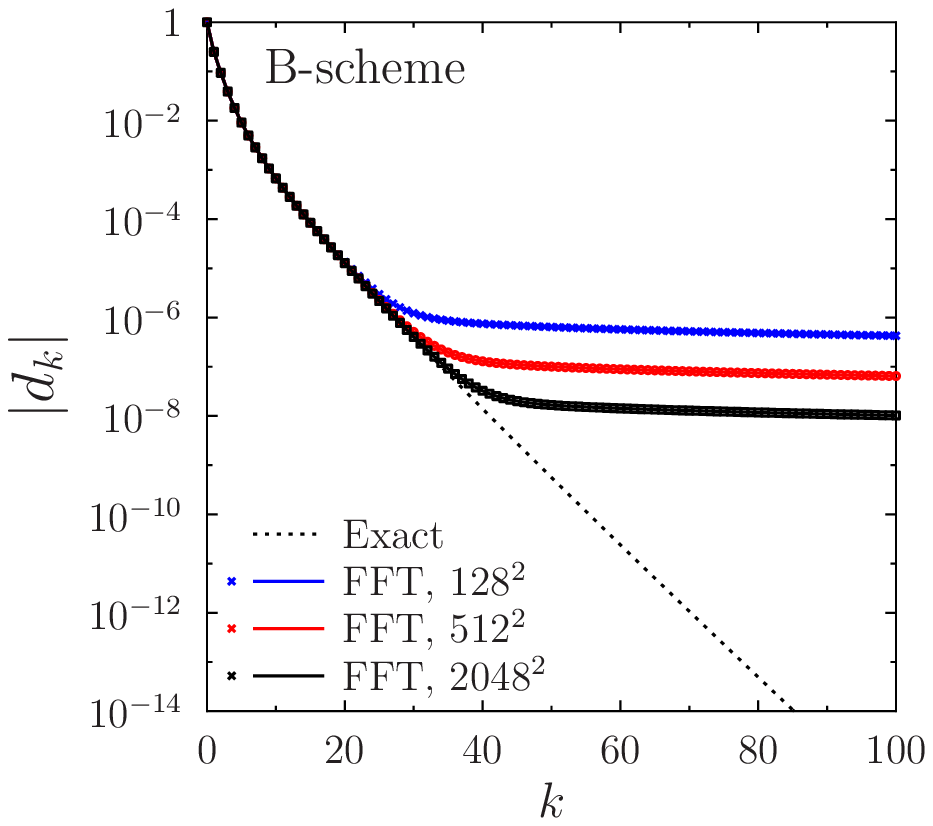} 
          \begin{center} \vskip -0.5cm  $\qquad$ (a)   \end{center}
    }
    \parbox{5.2cm}{
      \includegraphics[width=5.5cm]{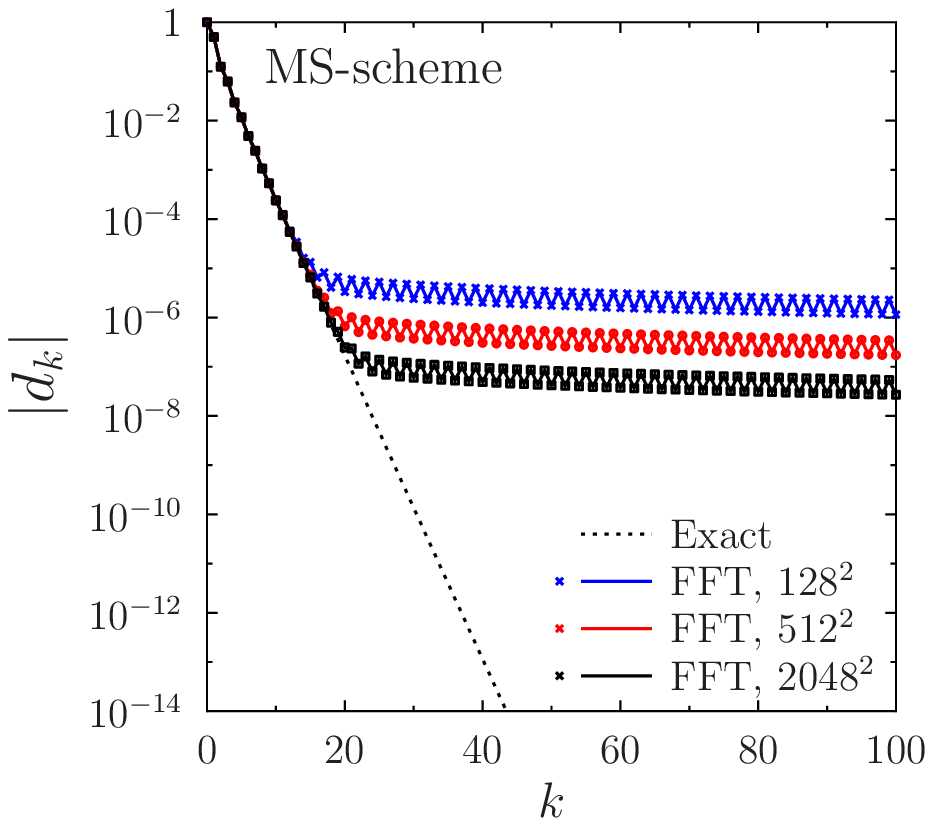}      
     \begin{center} \vskip -0.5cm $\qquad$ (b) \end{center}
    }
        \parbox{5.2cm}{
      \includegraphics[width=5.5cm]{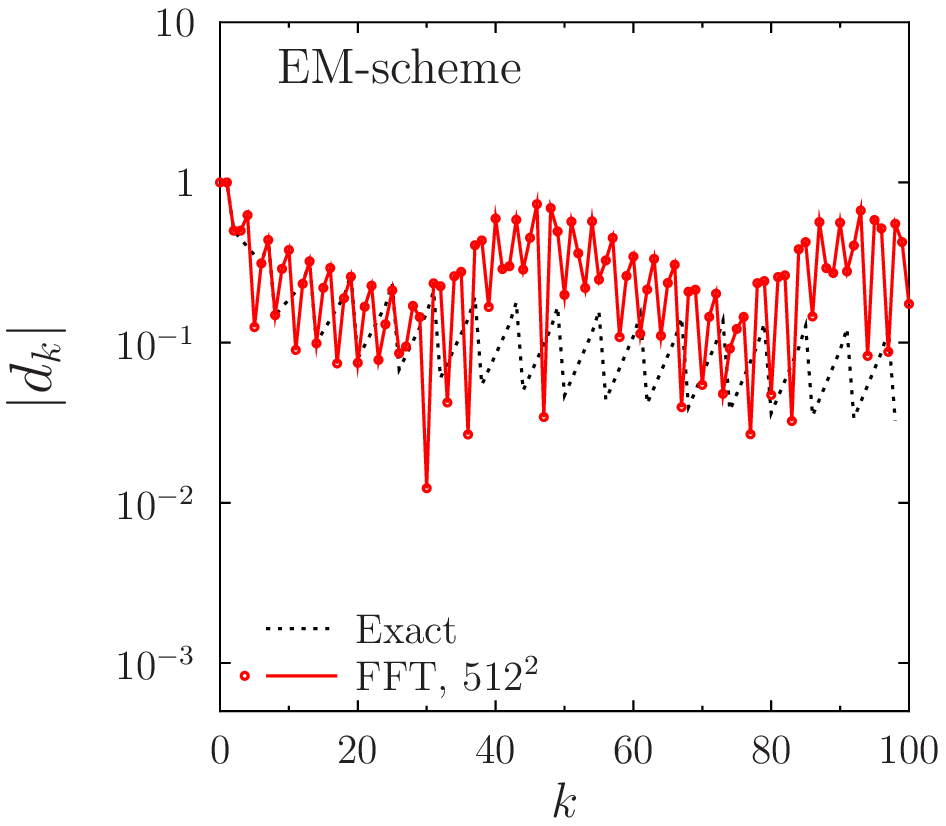}      
        \begin{center} \vskip -0.5cm $\qquad $ (c) \end{center}
    }
  \caption{
  Comparison between the coefficients {$d_k$} of the power series.
  a) Reference medium: matrix (B scheme). (b) Reference medium: arithmetic mean {(MS-scheme)}. (c) Reference medium: geometric mean ({EM-scheme}).
  }
   \label{comparaison3}
\end{figure}

\begin{figure}[!h]
  \begin{center}
  \parbox{5.cm}{
    \includegraphics[width=5.cm]{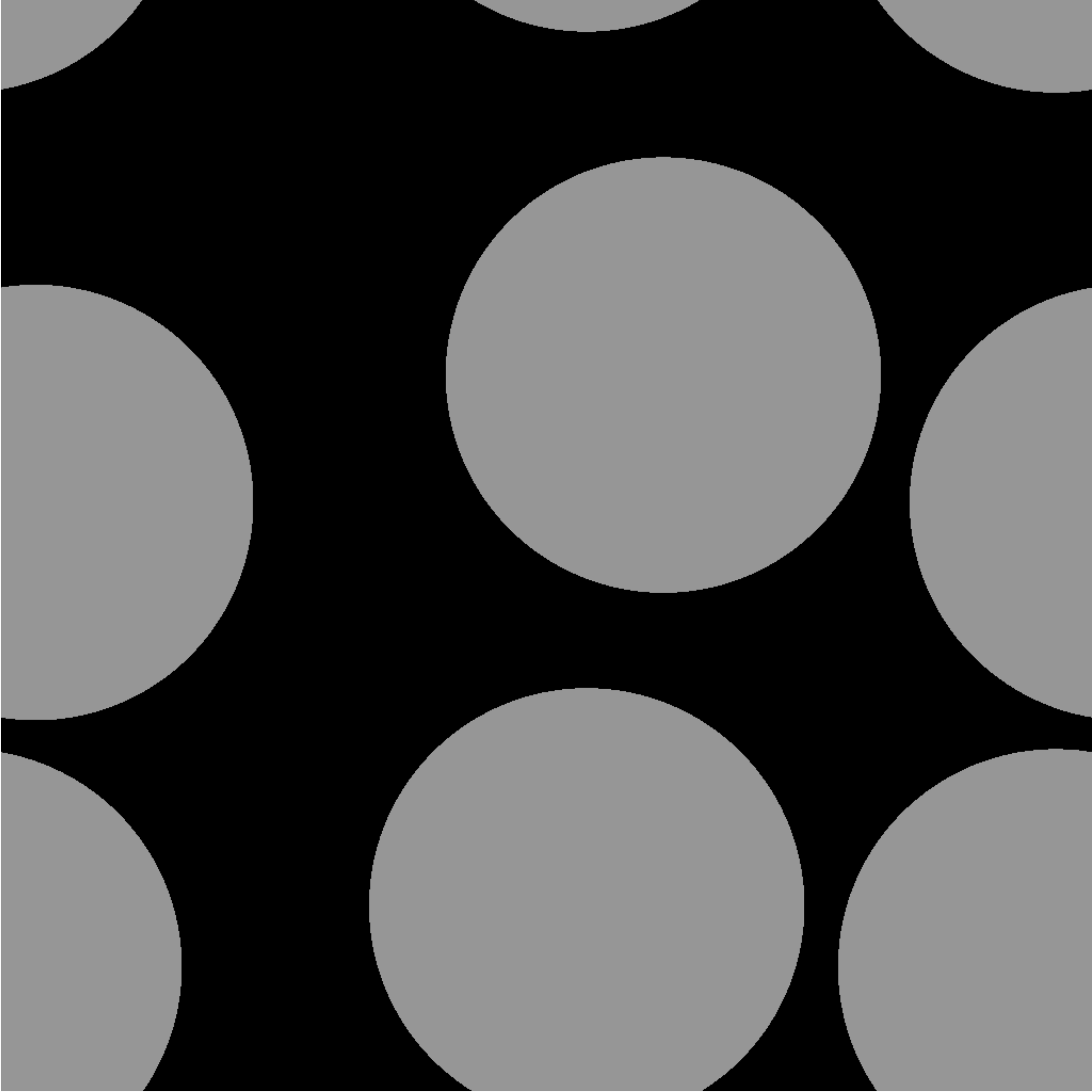}
  \begin{center} (a) \end{center}
      }
      \hskip 1cm
  \parbox{7cm}{
      \includegraphics[width=7cm]{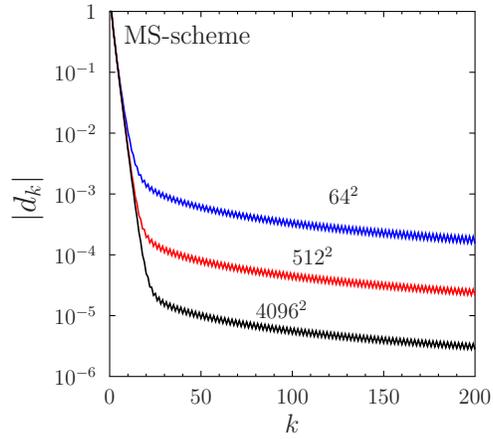}
       \begin{center} \vskip -0.5cm
      $\qquad$ (b) \end{center}
   }
    \end{center}
    \caption{ (a) Microstructure: 4 circular inclusions, volume fraction $50 \%$. (b) Coefficients  ${d}_k$ of the series {(\ref{recap2})} at different resolutions. MS scheme. 
    }
 \label{dk_4c}
\end{figure}

\subsection{Stopping criteria}\label{newconv}
\begin{figure}[!]
   \parbox{5.2cm}{
      \includegraphics[width=5.5cm]{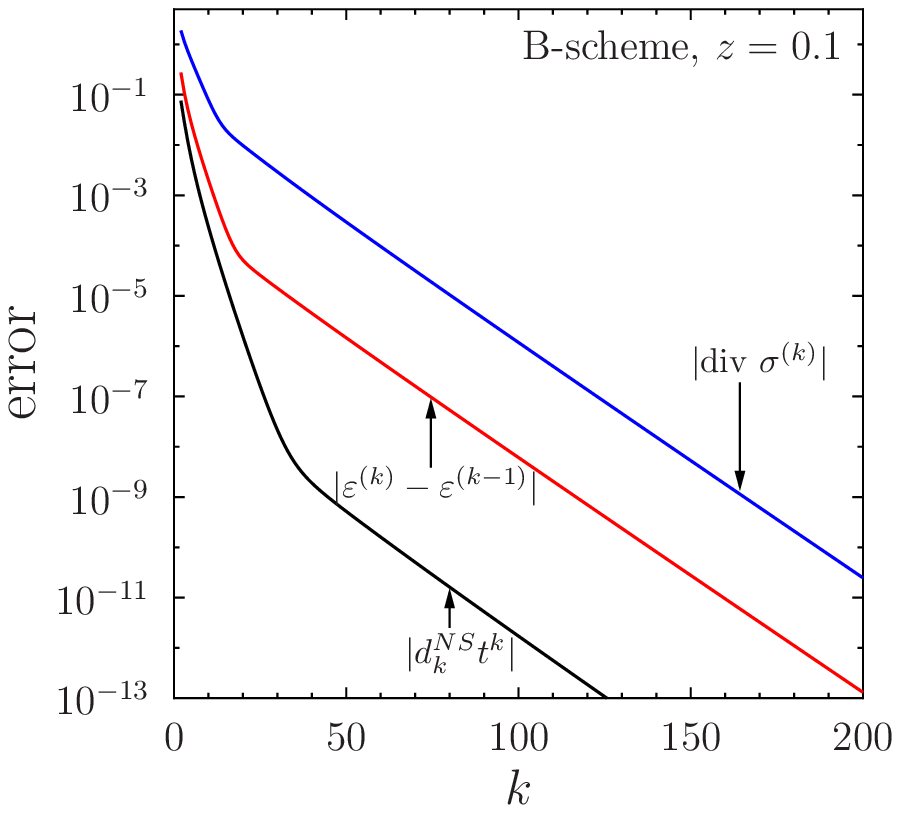}
      \begin{center} 
 $\qquad$     (a) \end{center}
    }
    \parbox{5.2cm}{
      \includegraphics[width=5.5cm]{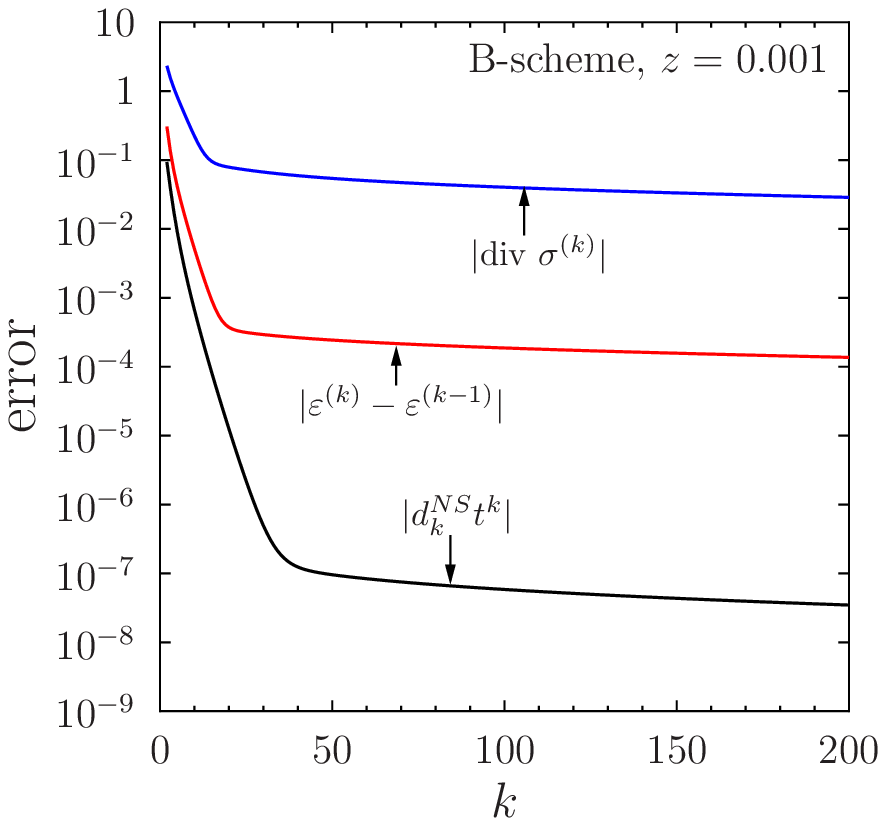}      
      \begin{center} 
$\qquad$      (b) \end{center}
    }
    \parbox{5.2cm}{
      \includegraphics[width=5.5cm]{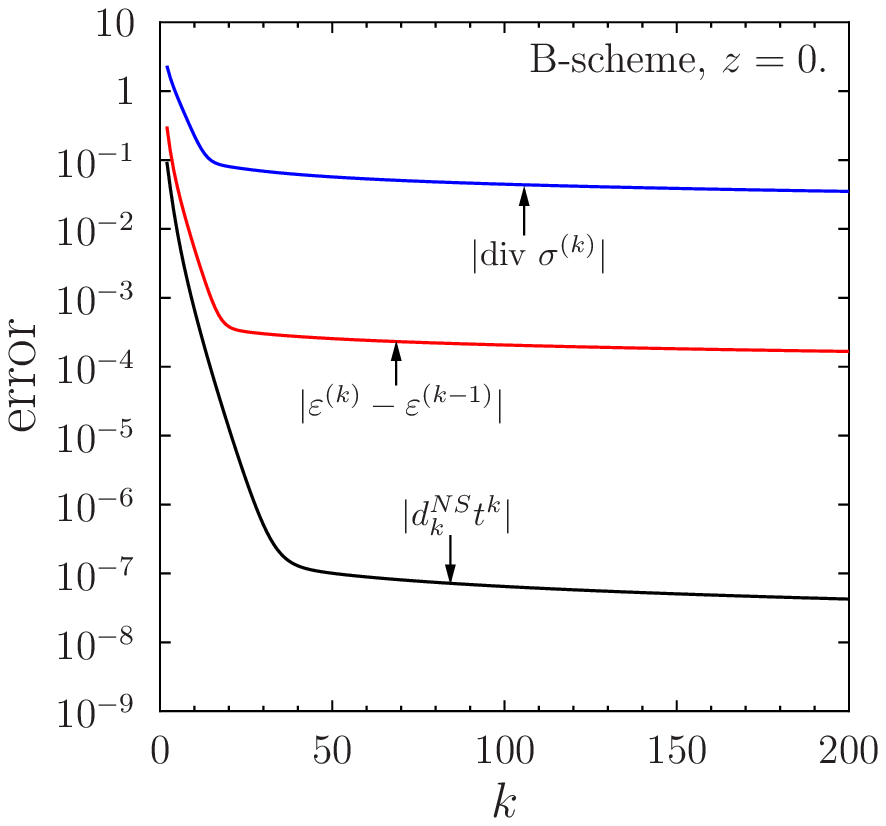}
      \begin{center} 
$\qquad$    (c) \end{center}
    }
  \caption{ Comparison between the error indicators (\ref{crit_equi}), (\ref {crit_diff}) and $|{d}^{NS}_k t^k|$. B-scheme (matrix as reference).  
   Three different contrasts: (a) $z=0.1$, (b) $z=0.001$ and (c) $z=0$.
    }
  \label{usual_error}
\end{figure}

The iterative procedure (\ref{LS4})$_1$ requires the choice of a criterion to stop the iterations, or alternatively to truncate the series (\ref{LS4})$_2$. Two error indicators which are popular in the literature are
 \begin{itemize}
  \item[-] equilibrium residual \citep{MOU98}: 
 \begin{equation}
\delta_1^{(k)} = \frac{ \norm{ \text{div}(\tod{\sig}^{(k)} )}_{L^2}}{\norm{\tod{\Sigma}}} = \frac{  \moy{ \abs{\text{div}(\tod{\sig} ^{(k)}(\tou{x}))}^2}^{1/2}} {\norm{\tod{\Sigma}}},
 \label{crit_equi}
 \end{equation}
where $\tod{\Sigma}$ is a stress which is natural in the problem at hand and used in (\ref{crit_equi})  to arrive at a non-dimensional criterion.
  The value commonly used for $\tod{\Sigma}$ is the average of the stress field: $\tod{\Sigma}=<\tod{\sigma}>$, but in the present study, for simplicity, it has been chosen with unit norm, i.e. $\norm{\tod{\Sigma}}=1$.
  \item[-] Difference between two  successive iterates \citep{KAB14}:
  \begin{equation}
 \delta_2^{(k)} = \norm{ \tod{\eps}^{(k+1)} - \tod{\eps}^{(k)}}_{L^2} = \moy { \abs{\tod{\eps}^{(k+1)} (\tou{x}) - \tod{\eps}^{(k)}(\tou{x}) }^2 }^{1/2}.
  \label{crit_diff}
  \end{equation}
  \end{itemize}
  Note that 
  $$\tod{\eps}^{(k+1)} - \tod{\eps}^{(k)}= -\toq{\Gamma}^0 \tod{\sig}^{(k)},$$ 
 from which it follows that the second criterion is equivalent to a criterion on the norm of $\toq{\Gamma}^0 \tod{\sig}^{(k)}$.
 \vskip 0.1cm
The iterative algorithm is stopped (or equivalently the Neumann series is truncated) as soon as the chosen error indicator is less than a given tolerance $\delta$:
$$ \delta_i^{(k)} \leq \delta, \; i=1 \; \text{or} \; i=2.$$ 
The variations with $k$ of these two error indicators are compared in figure \ref{usual_error} with the $k$-th term of the numerical power series  for the fixed-point scheme (\ref{LS4})$_1$ which can be used as a third error indicator
  \begin{equation}
  \abs{{d}_k  t^k}\norm{ \tod{E}}.
  \label{dketk}
\end{equation}
The following comments can be made. 
\begin{enumerate}
\item
 The most demanding criterion is the one based on equilibrium, {\em i.e.} the $L^2$ norm of $\text{div}(\tod{\sig})$. For a given error, it is the criterion requiring the larger number of iterations. Then comes the criterion based on the difference of two iterates (which is equivalently a check on equilibrium via the Green's operator $\toq{\Gamma}^0$). The third criterion \eqref{dketk}, based on the difference between two successive partial sums  in the Neumann expansion of the effective moduli is a global criterion, by contrast with the two others. It is therefore the less demanding criterion.
\item
The three curves for the three criteria are almost parallel. As observed in section \ref{discussion} the knee in the third curve (\ref{dketk}) is probably the point where the numerical error on the coefficients $d_k^{NS}$ due to discretization becomes significant and where the iterations should be stopped for this particular microstructure. It is interesting to note that the two other convergence criteria show the same knee. It is therefore expected that, whatever the criterion, all iterations beyond the knee do not improve the accuracy of the fields or of the effective properties. It is also expected that the threshold $\mathcal{K}$ would be of the same order for all three criteria. 
\item
A possible way to detect the threshold $\mathcal{K}$ , for any of these criteria, would be to conduct the simulations with two different discretizations and
 to stop iterating when the criterion for the two different discretizations deviate from each other. 
\end{enumerate}

\subsection{Discrete Green's operators}
\begin{figure}[htp!]
  \begin{center}
    \parbox{6.5cm}{
      \includegraphics[width=6.5cm]{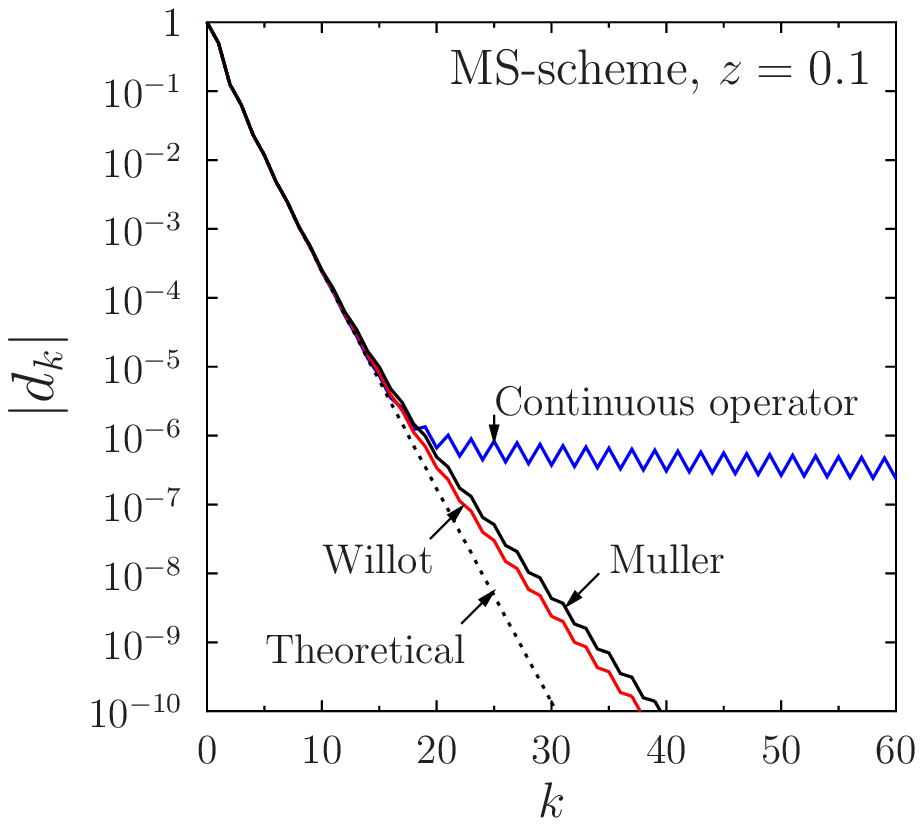}
        \begin{center} \vskip -0.5cm     $\qquad$(a) \end{center}
    }
    \parbox{6.5cm}{
      \includegraphics[width=6.5cm]{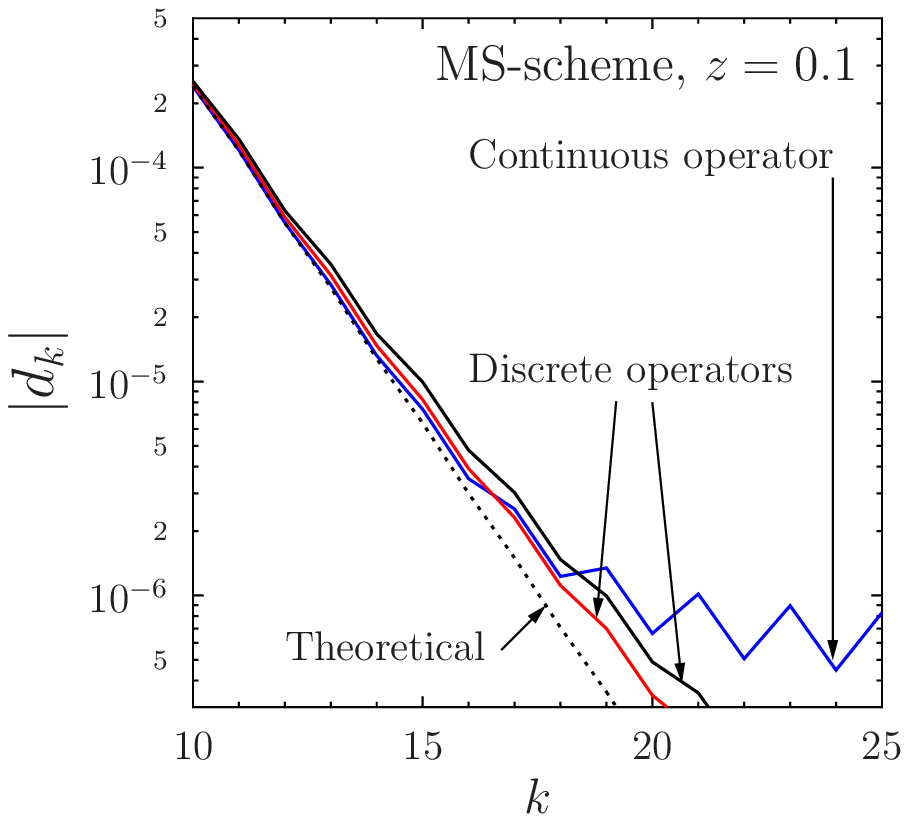}
 \begin{center} \vskip -0.5cm $\qquad$ (b) \end{center}
    }
  \end{center}  
  \vskip -0.5cm
    \caption{Obnosov's microstructure, discretization $512 \times 512$ pixels.   MS scheme.  Contrast $z=0.1$.
      Coefficients ${d}_k$ of the series {(\ref{series_scalar})} with different Green's operators: continuous operator (blue),   M\"uller's operator  \cite{MUL96} (black), {Willot-Pellegrini} operator \cite{WIL14} (red).  (a): iterations 1 to 60. (b): close-up on iterations 10 to 25.  
    }
  \label{willot1}
\end{figure}

Spurious oscillations are sometimes observed in the local fields when the continuous Green's operator {(\ref{Gamma0})} is used. Although certain authors attribute these oscillations to a  Gibb's phenomenon, they are more likely due to the fact that the derivative of a {\em nonsmooth function} (not $C^1$) is poorly approximated by the inverse Fourier transform of the Fourier transform of the function multiplied by the frequency. This has motivated several authors to introduce discrete Green's operators having a better behaviour with respect to derivation. 
Instead of using directly the expression of the continuous Green's operator (\ref{Gamma0})), the partial derivatives are approximated by finite differences
and incorporated into the Green's operator. Several variants have been proposed depending on the operator used to approximate the derivation (\cite{MUL96}, \cite{BRO02}, \cite{NEU02}, \cite{BRI10}, \cite{WIL14}, \cite{WIL15}, \cite{SCH16}). 
\vskip 0.1cm
The variations with $k$ of the coefficients ${d}_k$ are compared  in figure \ref{willot1} for three different choices of the Green's operator:  the continuous Green's operator used in Moulinec and Suquet \cite{MOU94, MOU98}, the modified Green's operators proposed by M\"uller \cite{MUL96} and by Willot and Pellegrini \cite{Willot:2008}. 

{
  The expressions used for the modified Green's operators are obtained by simply replacing in relation (\ref{Gamma0})
  the angular frequency $\tou{\xi}$ by $\tou{\xi}^M$ (as proposed by M\"uller), or $\tou{\xi}^W$ (as proposed by Willot  and Pellegrini),
  where
  {
    \begin{equation}
    \xi_j^{M} = \frac{N_j}{L_j} \sin \big( \frac{L_j}{N_j} \xi_j \big) \ 
    \quad
    (j=1,2 \ or \ 3),\quad 
    \xi_j^{W} = {\imath} \frac{N_j} {L_j} \Big( \exp \big( {-{\imath} \frac{L_j}{N_j} \xi_j} \big) - 1 \Big)
    \quad
    (j=1,2 \ or \ 3),
    \label{modified_xi}
    \end{equation}
  }
  (where $\xi_j$ is the j-th component of $\tou{\xi}$, ${\imath } = \sqrt{-1}$ and  $L_j$ and  $N_j$ are, respectively, the size and the number of pixels of the unit-cell in the $j$-th direction).
}

\begin{figure}[htp!]
\begin{center}
    \parbox{6.5cm}{
      \includegraphics[width=6.5cm]{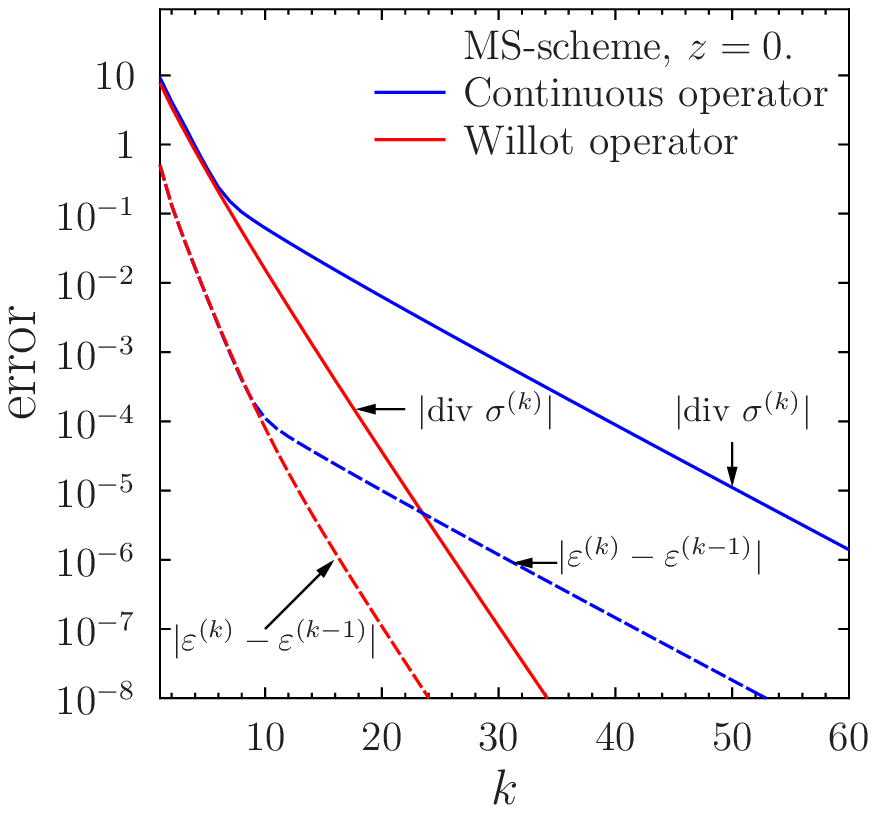}
      \begin{center}\vskip -0.5cm $\qquad$  (a) \end{center}
    }
    \parbox{6.5cm}{
      \includegraphics[width=6.5cm]{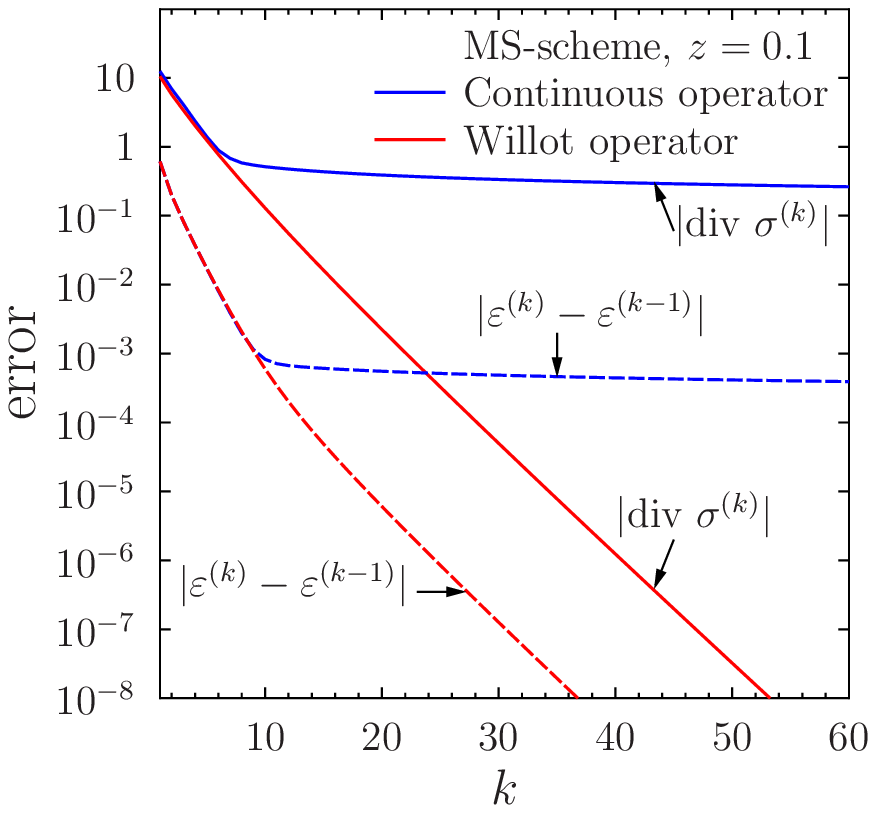}
      \begin{center} \vskip -0.5cm $\qquad$  (b) \end{center}
    }
\end{center}
\vskip -0.5cm
    \captionsetup{justification=raggedright}
    \caption{ Obnosov's microstructure, discretization $512 \times 512$ pixels.  MS scheme.  Error indicators (\ref{crit_equi}) (solid line), and (\ref{crit_diff}) (dashed line). Comparison between the continuous Green's operator (blue) and {Willot-Pellegrini}  discrete operator \cite{WIL14} (red).   (a):   Contrast $z=0.1$. (b) Contrast $z=0$.         }
  \label{willot3}
\end{figure}

\begin{figure}[!]
\begin{center}
    \parbox{6.5cm}{
      \includegraphics[width=6.5cm]{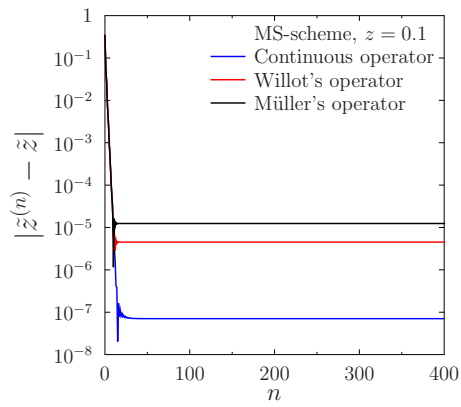}
      \begin{center} \vskip -0.5cm $\qquad$ (a) 
      \end{center}
    }
    \parbox{6.5cm}{
      \includegraphics[width=6.5cm]{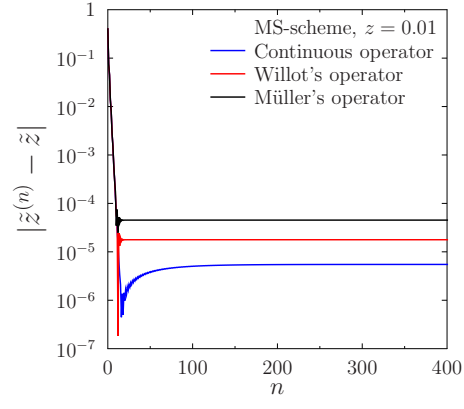}

      \begin{center} \vskip -0.5cm $\qquad$  (b) 
      \end{center}
    }

    \parbox{6.5cm}{
      \includegraphics[width=6.5cm]{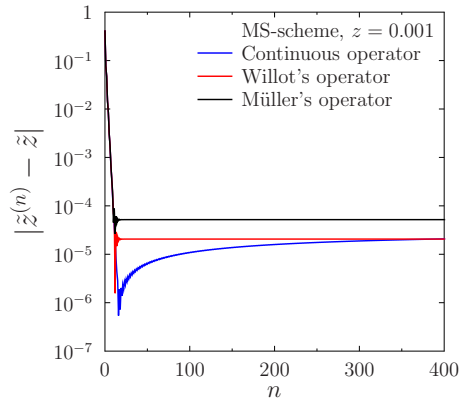}

      \begin{center}  \vskip -0.5cm $\qquad$ (c) 
      \end{center}
    }
    \parbox{6.5cm}{
      \includegraphics[width=6.5cm]{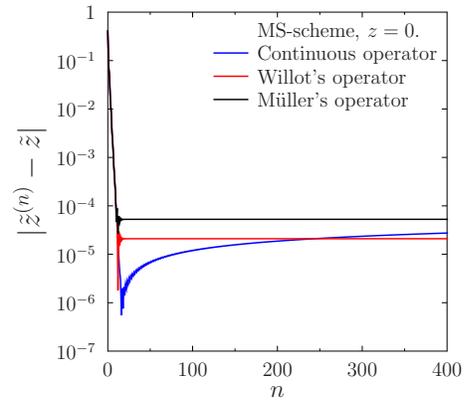}

      \begin{center} \vskip -0.5cm $\qquad$  (d) 
      \end{center}
    }
\end{center}
\vskip -0.5cm 
    \caption{  Obnosov's microstructure, discretization $512 \times 512$ pixels.  MS scheme.  Error on the effective property with different Green's operators: continuous operator (blue),  {Willot-Pellegrini}  operator \cite{WIL14} (red),  M\"uller's operator  \cite{MUL96} (black).  
(a):  $z=0.1$. (b): $z=0.01$. (c): $z=0.001$. (d): $z=0$.
  }
  \label{willot2}
\end{figure}

The MS scheme is used for this comparison but the same conclusions hold for the other schemes. The main observations are the following:
\begin{enumerate}
\item
In the first few iterations the predictions with the three Green's operator are in good agreement with {the} theoretical values, with slightly better predictions by the continuous Green's operator up to higher order (see the close-up in figure \ref{willot1}(b)). The continuous Green's operator seems to give the best estimation when $k$ is lower than 17, while the two discrete schemes deviate from the theoretical predictions at about $k=12$.
\item
  The numerical ${d}^{NS}_k$ deviate from the theoretical $d_k$ between iterations $k=15$ and $k=20$ for all three operators, the deviation occurring earlier for the discrete operators as already noticed. However, beyond that point, the predictions of the continuous Green's operator deviate more significantly from the theoretical coefficients and decrease only slowly, whereas the two modified Green's operators exhibit a more rapid decrease, although 
remaining much larger than it should be theoretically. As an illustration, in figure \ref{willot1} at iteration $k=30$, the coefficient obtained with Willot-Pellegrini  modified operator is approximately 18 times as large as the theoretical coefficient, the coefficient of Müller's operator is 32 times as large and the coefficient obtained with the continuous operator is 2772 times as large
as the theoretical coefficient.
In this sense the coefficients $d_k^{NS}$ for large $k$ are better predicted by the modified operators.
\item
This more accurate prediction of the $d_k$ beyond the threshold $\mathcal{K}$ goes together with a faster decrease of the two classical error indicators 
(\ref{crit_equi}, \ref{crit_diff}) for the discrete operators {as shown in {figure \ref{willot3}}. It is therefore very likely that the stopping criterion based on the discrete Green's operator (all iterates being computed with this discrete operator {and the equilibrium tests being performed using the modified
    frequencies  \eqref{modified_xi}}) is less sensitive to the discretization error than the continuous operator. }
\item
These observations seem to favour the use of discrete Green's operators. However, the comparison of the effective properties predicted by the different Green's operators in  figure \ref{willot2} goes in the opposite direction. The effective properties computed with the continuous Green's operator are in better agreement with the theoretical value, except in the extreme case $z=0$. This can be understood by going back to the expansion (\ref{Dk}). The first coefficients of the series are better predicted by the continuous operator up to a threshold $\mathcal{K}$ which is larger than for the discrete operators. And although the coefficients for larger $k$ are better predicted by the discrete operators, these coefficients enter the expansion (\ref{Dk}) with a factor $t^k$ which is very small for large $k$ when $\abs{t} <1$. So the coefficients which are crucial for an accurate prediction of the effective property are the coefficients of lower order, which are better predicted with the continuous operator. The gain in accuracy provided by the discrete operators for the coefficients of higher order is lost by the multiplication by $t^k$ when $\abs{t} <1$. However when $t$ approaches 1, or even become larger than 1, these coefficients become again crucial and the predictions of the discrete operators become more accurate. As a tentative conclusion, one could say that the continuous Green's operator performs better for moderate contrast whereas discrete operators should be used for large contrast.
\end{enumerate}

\newpage
\subsection{Theoretical rate of convergence of the three iterative schemes}\label{sec:rate}
\begin{figure}[htp!]
\begin{center}
    \parbox{5.cm}{
             \includegraphics[width=5cm]{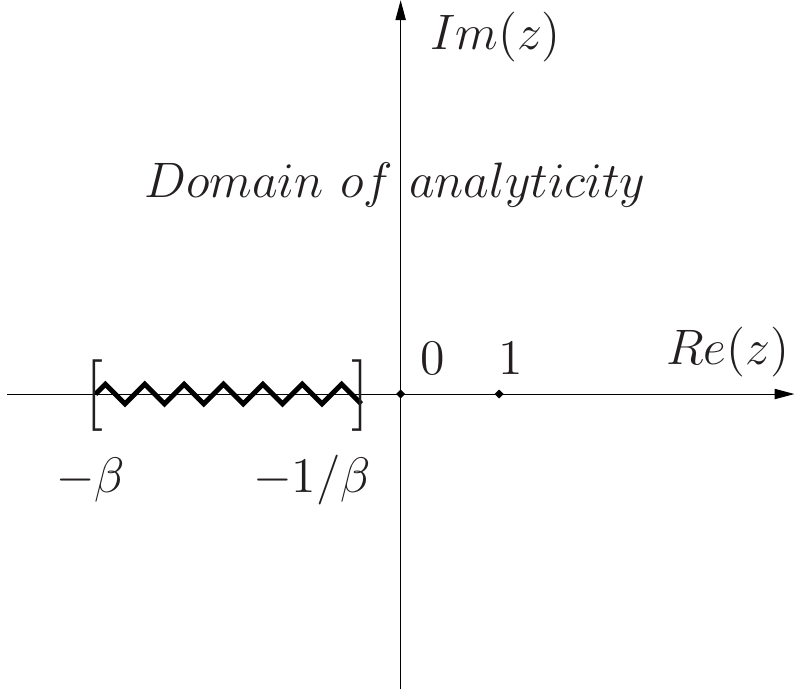}
    }
   \end{center} 
 \caption{Domain of analyticity of $\eff{z}$.}
  \label{definition}
\end{figure}

\begin{figure}[htp!]
    \begin{center}
    \parbox{5.cm}{
            \includegraphics[width=5cm]{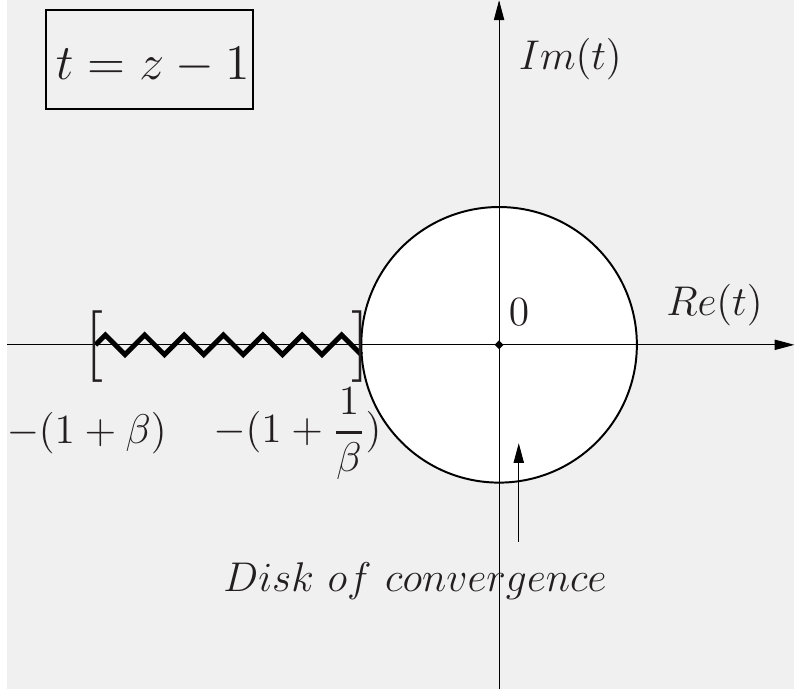}
          \begin{center} \vskip -0.5cm
      (a) $\quad \quad$ \end{center}
    }
    \hskip 0.3cm
    \parbox{5.cm}{
         \includegraphics[width=5cm]{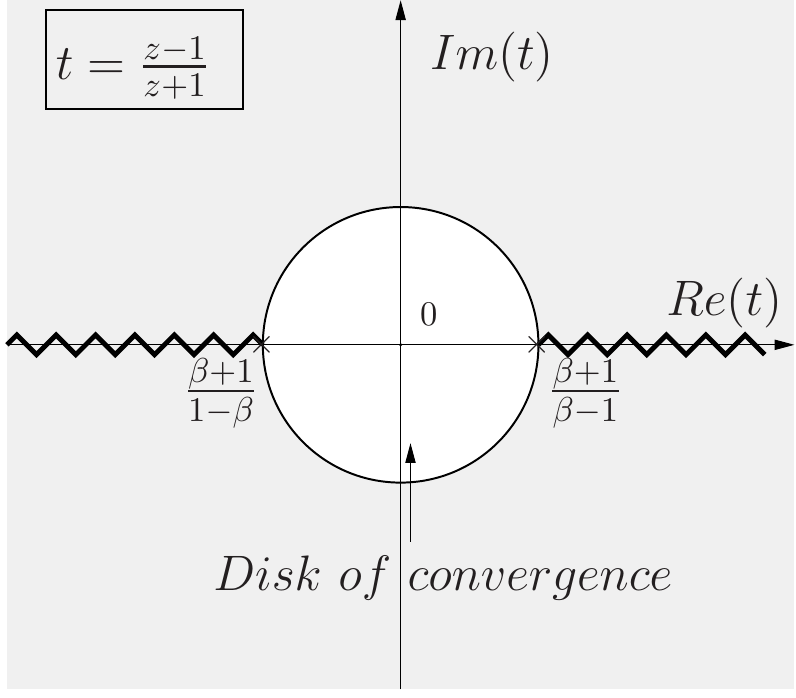}     
                   \begin{center}  \vskip -0.5cm   (b) \end{center}
    }
     \hskip 0.3cm
    \parbox{5.cm}{
        \includegraphics[width=5cm]{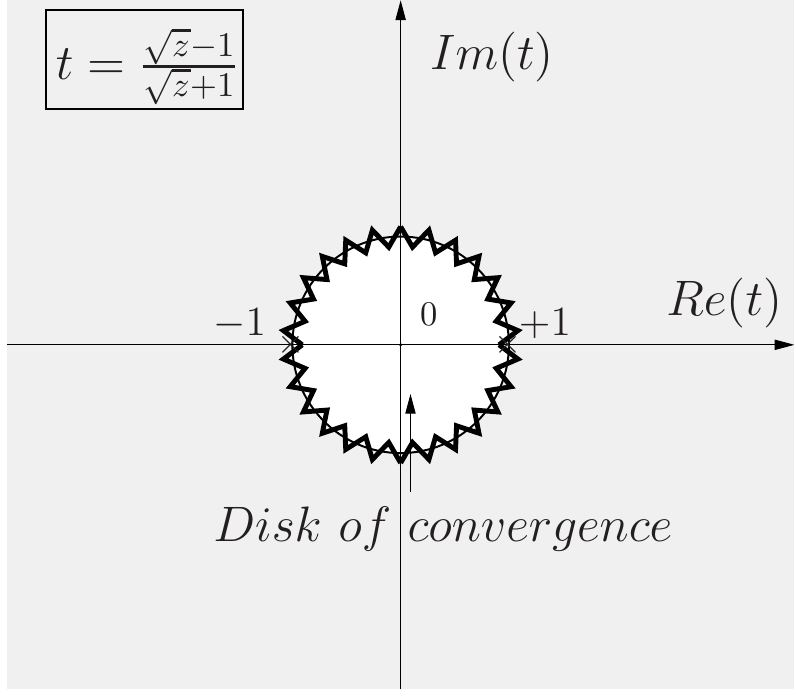}      
            \begin{center} \vskip -0.5cm   $\qquad \qquad$ (c) \end{center}
    }
    \end{center}
    \caption{Domain of convergence of the power series in the neighborhood of $t=0$. (a) $t=z-1$ (B-scheme. (b): MS-scheme. (c): EM-scheme. 
     }
  \label{domain}
\end{figure}

We end this section by a discussion of the (theoretical) convergence rate of the three iterative schemes (Brown, Moulinec-Suquet, Eyre-Milton) for microstructures which are more general that Obnosov's microstructure. 
These rates of convergence depend on the microstructure and on the contrast between the phases. The present section shows  that none of these iterative methods will converge faster than the others for all possible microstructures and phase contrasts. The choice of the fastest method depends on the material data. Therefore any comparison between these methods (and others), as often reported in the literature, should be taken with care, as the conclusions may depend on the microstructure under consideration and on the contrast between the phases. 
\vskip 0.1cm 
This dependence can be understood on a simple example, that of the effective conductivity of  a two-phase composite with two isotropic phases.  The singularities of the complex conductivity $\eff{z}$ as a function of $z$ are located on the real negative axis. It is further assumed that these singularities lie in an interval $[-\beta, -\frac{1}{\beta}]$ with $\beta \geq 1$ (see figure \ref{definition}). Specifically, we assume that there is a singularity at $z=-1/\beta$
and no singularities for $z<-\beta$ (implying no pole at $z=\infty$).
Obnosov's microstructure corresponds to $\beta=3$. The checkerboard microstructure corresponds to $\beta=+ \infty$.
\vskip 0.1cm 
It is to be cautioned that the analysis here is the theoretical rate of convergence in the asymptotic regime where the number of iterates is very large. In particular, suppose
the singularity at $z=-1/\beta$ is say a pole with extremely small residue, and there is a pole at $z=-1/\gamma$ with a significant residue, and no singularities
in $(-1/\gamma,-1/\beta)$, nor for { $z<-\gamma$ ($\gamma$ is supposed to be larger than 1)}. Then the initial rate of convergence of the iterates should be dictated by the pole at $z=-1/\gamma$, with the pole
at $z=-1/\beta$ only slowing down the convergence rate after many iterates (when the error is of the same magnitude as the effect due to the extremely small residue
of the pole at $z=-1/\beta$). One then expects a knee in the convergence rates, similar to our numerical simulations. As in our numerical simulations, such knees 
can also be due to spurious singularities caused by the discretization, and these could swamp the effect of an extremely small residue at $z=-1/\beta$. 
We ignore such considerations in the following analysis, thus
assuming the singularity at $z=-1/\beta$ to be significant enough to control the convergence rate during the crucial initial stage of iterations.
The three schemes  (Brown, Moulinec-Suquet, Eyre-Milton) and the corresponding iterative Neumann series are associated with a contrast variable $t$ which reads in the three respective cases
\begin{equation}
 \text{Brown:}\; t=z-1,\quad \text{Moulinec-Suquet:}\;  t= \frac{z-1}{z+1},\quad  \text{Eyre-Milton:}\; t= \frac{\sqrt{z}-1}{\sqrt{z}+1}.
 \label{transform}
 \end{equation}
 The domains of  $z$ for which the  iterative schemes converge is the radius of convergence of the power series in the $t$-plane in the neighborhood of $t=0$. It is the largest disk contained in the domain of analyticity of $\eff{z}$ as a function of $t$. These domains, schematically represented in figure \ref{domain}, are respectively
\begin{equation}
 \mathcal{D}_B = \{ \abs{t} < 1 + \frac{1}{\beta} \}, \quad 
  \mathcal{D}_{MS} = \{ \abs{t} < \frac{\beta+1}{\beta-1} \},\quad   \mathcal{D}_{EM} = \{ \abs{t} < 1  \}.
 \label{domains}
 \end{equation}
The radius of convergence of the corresponding series in power of $t$ are given 
\begin{equation}
\rho_B= 1 + \frac{1}{\beta}, \quad \rho_{MS}=\frac{\beta+1}{\beta-1}, \quad \rho_{EM}=1.
\label{rate1}
\end{equation}
The three power series $\displaystyle \sum_{k=0}^{+\infty} d_k t^k$  (and the associated iterative schemes) converge when $\abs{t} < \rho$, {where}  $ t=z-1,\, \frac{z-1}{z+1},\; \text{and}\;  \frac{\sqrt{z}-1}{\sqrt{z}+1}$, and $ \rho=\rho_B,\,  \rho_{MS},\, \rho_{EM}$ {respectively}. It is therefore natural to evaluate, for a given contrast (characterized by $z$) and a given microstructure (associated with $\beta$), 
 the rate of convergence of a series expansion through the ratio $r = \rho/\abs{t}$. The larger $r$, the faster the convergence of the series expansion (and of the associated iterative schemes).  This ratio reads, for each choice of the variable $t$:
\begin{equation}
r_B= \abs{\frac{\beta+1}{\beta(z-1)}}, \quad r_{MS}=\abs{\frac{(\beta+1)(z+1)}{(\beta-1)(z-1)}}, \quad r_{RM}=\abs{\frac{\sqrt{z}+1}{\sqrt{z}-1}}.
\label{rate}
\end{equation}
\begin{figure}[!]
\includegraphics[width=18cm]{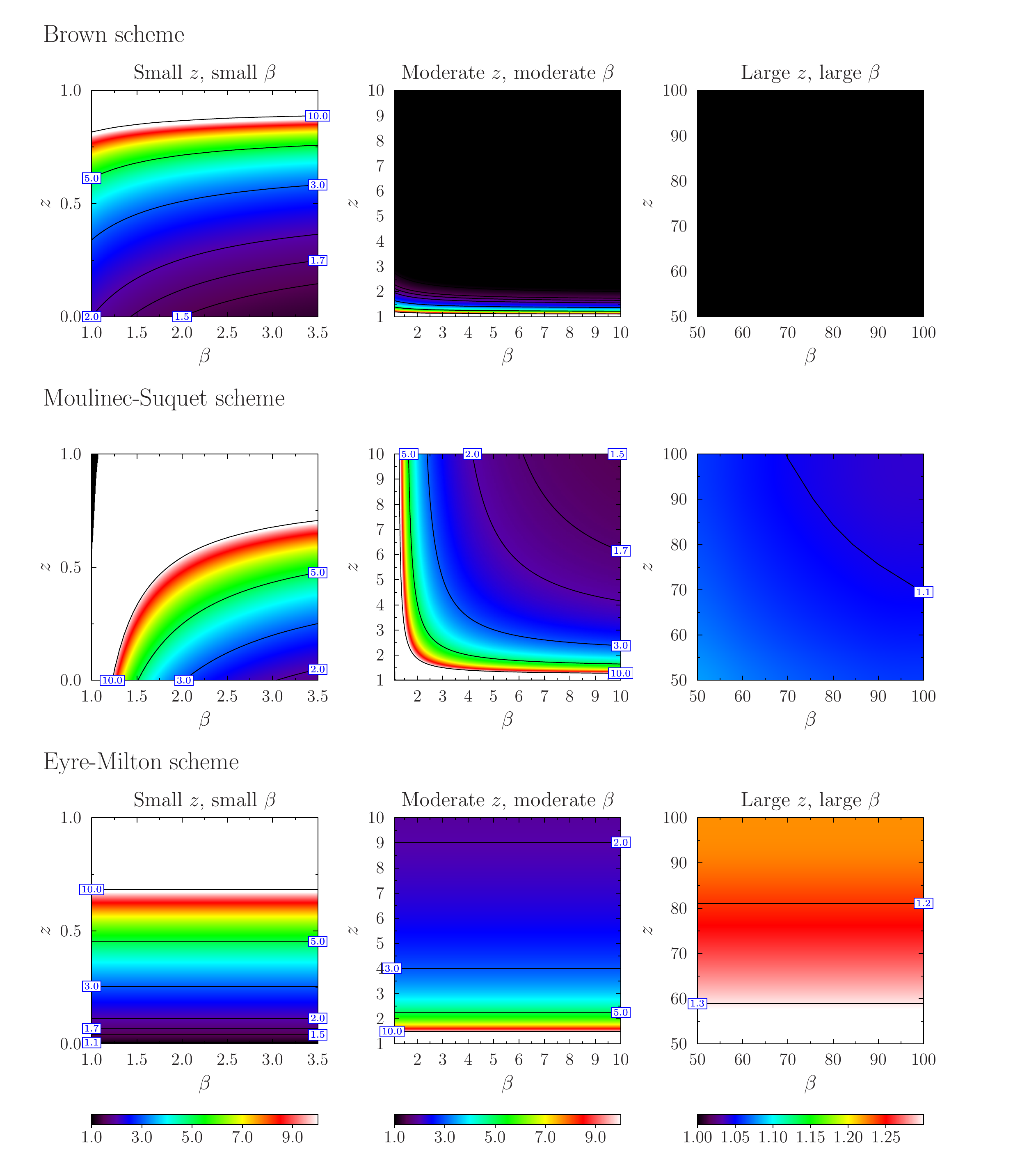}
\caption{Snapshots of the rate of convergence for the 3 schemes in the plane $(\beta,z)$. The brighter the color, the faster the rate of convergence.}
\label{fig3}
\end{figure}
Snapshots of the three ratios are shown in Figure \ref{fig3}, for small, moderate and large values of $z$ and $\beta$. The snapshots corresponding to a given scheme are placed horizontally. The three schemes can be compared by looking vertically at the columns, the color table at the bottom being the same for a given column. The larger the ratio, the better the rate of convergence.
\vskip 0.1cm
The ratios can be compared analytically two-by-two in the domains of $z$ where they are defined. Attention will be limited to real values of $z$. 
\paragraph{Brown vs Moulinec-Suquet}
The joint domain of definition  of the two power expansions in the $(\beta,z)$-plane is 
\begin{equation}
z > - \frac{1}{\beta}, \quad \beta \geq 1.
\label{domainMSB}
\end{equation}
Straightforward algebra shows that 
\begin{equation}
\frac{r_{MS}}{r_{B}} = \abs{ \frac{\beta z+1}{\beta -1}+1 } > 1,
\label{rate_MSB}
\end{equation}
where the last inequality comes from (\ref{domainMSB}). In other words, the MS-scheme always over-performs the B-scheme. 
\paragraph{Brown vs Eyre-Milton}
The joint domain of definition of the two power expansions in the $(\beta,z)$-plane is 
\begin{equation}
z \geq 0, \quad \beta \geq 1.
\label{domainEMB}
\end{equation}
After straightforward algebra, one obtains that  
\begin{equation}
\frac{r_{EM}}{r_{B}} = (\sqrt{z}+1)^2 \frac{\beta}{\beta+1},
\label{rate_EMB}
\end{equation}
and
\begin{equation}
r_{B} \geq  r_{EM} \quad\text{when}\; 0 \leq z \leq z_0= \left(\sqrt{\frac{\beta +1}{\beta}}  -1\right)^2, \quad r_{EM} \geq  r_{B} \quad\text{when}\; z \geq z_0. 
\label{rate_EMB2}
\end{equation}
In other words, for small enough $z$ (less than $z_0$) the Brown scheme over-performs the Eyre-Milton scheme, but the Eyre-Milton scheme is faster when $z \geq z_0$. For the Obnosov's microstructure $\beta=3, z_0 \simeq 0.23932$. 
\begin{figure}[!]
\begin{center}
    \parbox{6.cm}{
             \includegraphics[width=6cm]{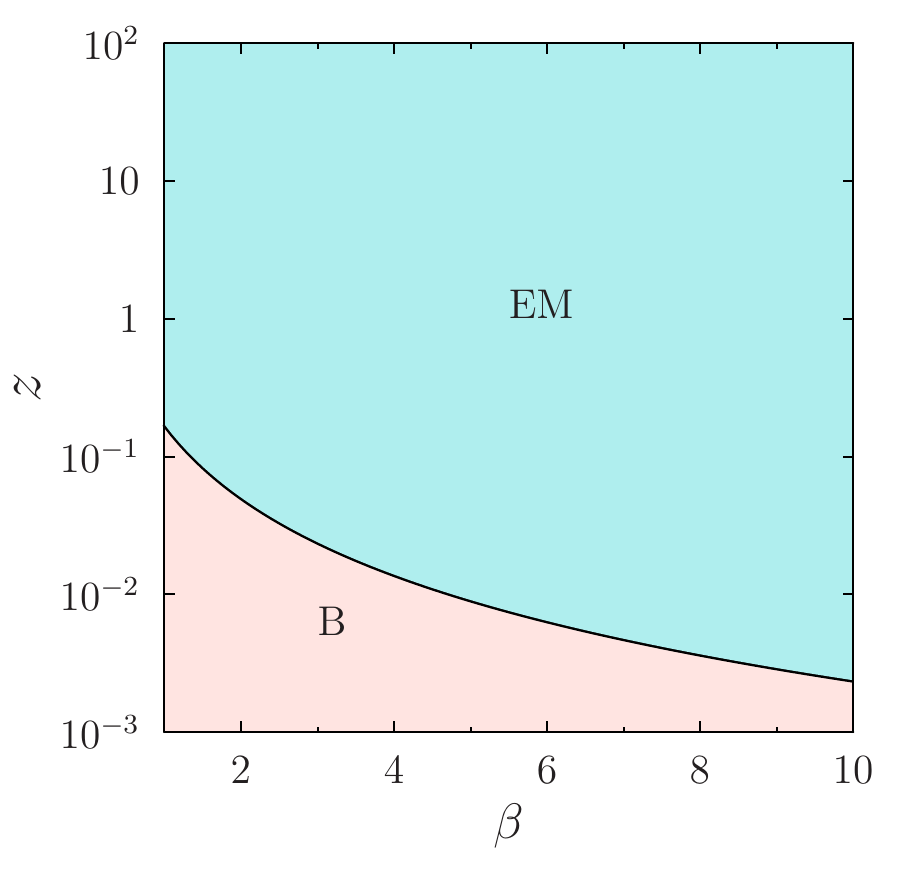}
          \begin{center} \vskip -0.5cm
      (a) $\quad \quad$ \end{center}
    }
    \hskip 0.3cm
    \parbox{6.cm}{
          \includegraphics[width=6cm]{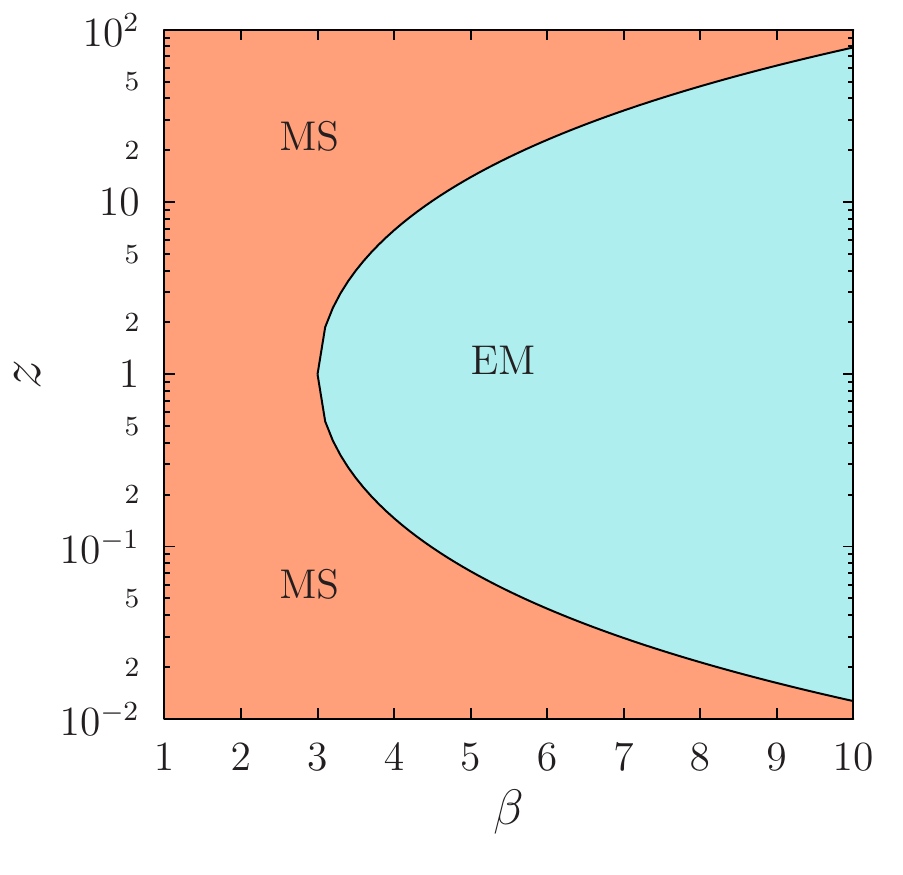}     
                   \begin{center}  \vskip -0.5cm   (b) \end{center}
    }
    \end{center}
    \caption{Best rate of convergence. (a) EM-scheme vs B-scheme. (b): All 3 schemes.}
  \label{3schemes}
\end{figure}
\paragraph{Eyre-Milton vs Moulinec-Suquet}
The joint domain of definition  of the two power expansions in the $(\beta,z)$-plane is 
\begin{equation}
z \geq 0, \quad \beta \geq 1.
\label{domainMSEM}
\end{equation}
Straightforward algebra shows that 
\begin{equation}
\frac{r_{EM}}{r_{MS}} = \displaystyle \frac{ 1  + \frac{2 \sqrt{z}}{z+1}}{ 1 + \frac{2}{\beta -1}},
\label{rate_EMMSa}
\end{equation}
Two different cases are found
\begin{enumerate}
\item
When  $1 \leq \beta \leq  3$,  then 
$$ \frac{2}{\beta -1} \geq 1 \geq \frac{2 \sqrt{z}}{z+1}\quad \Rightarrow\quad 
\frac{r_{EM}}{r_{MS}}  \leq 1,$$ and the MS-scheme is faster than the EM-scheme, regardless of $z$. The only case where ${r_{EM}}/{r_{MS}}  = 1$ is when $\beta=3$ and $z=1$ (homogeneous material).
\item
When $\beta >  3$ then
\begin{eqnarray}
r_{EM} \geq  r_{MS} \quad\text{when}\; z_1 \leq z \leq z_2, \quad r_{MS} \geq  r_{EM} \quad\text{when}\; 0 \leq z \leq z_1 \; \text{or}\; z \geq z_2, \nonumber \\[2ex]
\text{with}\; z_1= \frac{1}{4} \left[\beta-1 - \sqrt{(\beta-1)^2-4} \right]^2,\quad z_2= \frac{1}{4} \left[\beta-1 + \sqrt{(\beta-1)^2-4} \right]^2.
\label{rate_EMMSb}
\end{eqnarray}
\end{enumerate}
The scheme offering the best ratio is therefore dependent on the problem at hand. It depends on the microstructure through $\beta$ and on the contrast between the phases through $z$. A map of the comparison between the 3 schemes is given in figure \ref{3schemes}. 
When $\beta$ is large, then $z_1 \sim 1/(\beta-1)^2$ tends to 0, and $z_2 \sim (\beta-1)^2$ tends to $+\infty$. In particular when $\beta=+\infty$ (checkerboard for instance), the EM-scheme has a better ratio of convergence than the MS-scheme. This confirms the column on the right of figure \ref{fig3}. 
\vskip 0.1cm
{{\bf Remark 4:} As we have just seen, the
relative performance of the different methods depends on the conductivity ratio $z$ and on the interval on the negative real $z$-axis 
outside of which there are no singularities in the conductivity function (for the Obnosov microstructure this is the interval $[-3,-1/3]$). If this interval is known 
the new accelerated FFT algorithm described in Milton \cite{Milton:2016:ETC} chap. 8, should be superior to all three methods. However it is to be noted that calculating, or estimating, this interval for an arbitrary geometry is a highly nontrivial problem. Some progress was made in Bruno \cite{Bruno:1991:ECS}.} 

\section{Conclusion}
The starting point of this study is the {classical} observation that the convergence of certain iterative methods used to invert the Lippmann-Schwinger operator in heterogeneous media, is directly related to the convergence of series expansions of the effective properties in powers of a contrast variable. Three different such Neumann series are  investigated in details, first from a theoretical perspective  and then numerically. A specific example (Obnosov's microstructure) where effective properties are available in closed form serves for the comparison. 
\vskip 0.1cm
A first result {of the present study is that the range of convergence of the series given in the literature with no information on the microstructure (recalled here in section 2) can be extended when additional information about the microstructure is available (section 3). In particular, when the location of the singularities of the effective moduli in the complex plane of phase contrast is known, the theoretical prediction of the range of contrast for which the iterative methods converge and the prediction of the rate of convergence of these methods can be improved significantly. Unfortunately these theoretical improvements are not always confirmed by the numerics.} 
\vskip 0.1cm
Second, {in an attempt to understand  this discrepancy between the theory and the numerics,} it is noted that the approximation of the field solutions in a finite dimensional space introduces a systematic error on the effective properties (because of their variational character). A direct consequence of this error is that beyond a certain order $\mathcal{K}$ the numerical coefficients of the series differ from their exact value. This is confirmed by the observation that the deviation {between the theory and the numerics} occurs later when the discretization is refined.  The theoretical results for convergence are for a continuous problem, whereas the numerical results are for a discrete, finite-dimensional approximation of it. The error on the coefficients of high order in the series do not influence much the convergence of the series at moderate contrast, but has a dramatic influence at high contrast. This explains why Neumann series fail to converge numerically for certain values of the contrast for which they should converge theoretically.
\vskip 0.1cm
{Thirdly, The rate of convergence of the three iterative schemes are compared when the contrast $z$ between the phases and a parameter $\beta$ related to the location of the singularities of the effective moduli are varied. It is found that none of the three schemes over-performs the others for all values of $z$ and $\beta$. The fastest scheme depends on the microstructure and on the phase contrast.} 
 \vskip 0.1cm
{\ajj{Fourth}, when discrete Green's operators are used, the detection of the moment where the iterations should be stopped seems to be more accurate. This does not always result in an improvement of the effective properties but it is likely that the local} {fields} {are better captured by avoiding superfluous iterations which may only add noise to the fields. }
\vskip 0.1cm
\ajj{We end by remarking that while we have restricted our attention to two-phase composites, for composites with more than two phases the convergence of the various approximation schemes is again
dictated by the singularities of the effective conductivity as a function of the component conductivities. The analytic properties of such functions have been studied in
\cite{GPMult} and \cite{Milton:1990:RCF}.}

\ajj{\subsection{Acknowledgements}
Fruitful discussions with C. Bellis and F. Silva are gratefully acknowledged. HM and PS have received funding from Excellence Initiative of Aix-Marseille University - A*MIDEX, a French "Investissements d'Avenir" programme  in the framework of the Labex MEC. {GWM thanks the National Science Foundation for support through  Grant DMS-1211359}.
}

\newpage
\appendix
{\section{The Green's operator $\toq{\Gamma}^0$ and useful properties of the key operator $\toq{\Gamma}^0\toq{L}^0$}}\label{appendixA}
The space of periodic $L^2$ tensor fields on $V$
 is denoted by $\mathscr{H} = \tou{\mathcal{L}}^2_{\sharp}\left(V\right)$. 
  It is a Hilbert space when equipped with the scalar product
$$
 \moy{\tod{e},\tilde{\tod{e}}}  = \frac{1}{\abs{V}} \int_V \tod{e}(\tou{x}). \tilde{\tod{e}}(\tou{x}) d\tou{x},  
 $$
 with 
 $$ \tod{e}. \tilde{\tod{e}}= \sum_{i=1}^N e_{i} \tilde{e}_{i}\; \text{for 1st order tensors, or}\quad \tod{e}. \tilde{\tod{e}}= \sum_{i,j=1}^N e_{ij} \tilde{e}_{ij} \; \text{for 2nd order tensors}.$$
{
Let $\toq{L}^0$ denote a uniform (no spatial dependence) fourth-order tensor, positive definite with major and minor symmetries. For a given field $\tod{\tau}$ in  $\mathscr{H}$, consider the following {Eshelby} problem:
\begin{equation}
\textit{Find } \tod{\sig}\in  \mathscr{S} \textit{ and } \tod{\eps}^*\in  \mathscr{E}_0 \textit{ such that: }
 \tod{\sig} = \toq{L}^0: \tod{\eps}^*-\tod{\tau}, 
 \label{thermoelas}
\end{equation}
where
\begin{equation}
\left.
\begin{array}{l}
\mathscr{E}_0 = \Big\{ \tod{\eps}^*\in \mathscr{H}  \text{ such that:}\  \exists\, \tou{u}^* \in \tou{H}^1_{\sharp}\left(V\right) ,\ \tod{\eps}^* = \frac{1}{2} (\tod{\nabla u}^* + \tod{\nabla} {\tou{u}^*}^\top) \Big\}, \\[2ex] \mathscr{S} = \Big\{ \tod{\sig} \in \mathscr{H} \text{ such that:}\  \hbox{\rm div }\tod{\sigma}(\tou{x})=\tou{0} \text{ in }\ V, \  \tod{\sig}.\tou{n} \text{ anti-periodic on }\ \partial V\Big\}.
\end{array}
\right\}
\label{perio}
\end{equation}
The problem (\ref{thermoelas}) has a unique solution $\tod{\sig} \in \mathscr{S}$ and $\tod{\eps} \in \mathscr{E}_0$ and the periodic strain Green's operator $\toq{\Gamma}^0$ of the reference medium with stiffness $\toq{L}^0$ is defined as 
$$
\toq{\Gamma}^0:\   \tod{\tau} \in  \mathscr{H} \longmapsto\toq{\Gamma}^0\tod{\tau} = \tod{\eps}^* \text{ solution of (\ref{thermoelas})}.
$$}
$\mathscr{H}$ can be alternatively equipped with an energetic scalar product 
\begin{equation}
\mmoy{ \tod{e},\tilde{\tod{e}}}= \moy{\tod{e}:\toq{L}^0:\tilde{\tod{e}}} , 
\label{scalar}
\end{equation}
\noindent
 The operator $\toq{\Gamma}^0 \toq{L}^0$ has the following  properties on $\mathscr{H}$ endowed with the energetic scalar product (\ref{scalar}).
\begin{enumerate}
\item $\toq{\Gamma}^0 \toq{L}^0$ is a self-adjoint operator from $\mathscr{H}$ endowed with the scalar product (\ref{scalar}) into itself:
\begin{equation}
\mmoy{\tod{e}_1, \toq{\Gamma}^0 \toq{L}^0 \tod{e}_2} = \mmoy{ \toq{\Gamma}^0 \toq{L}^0\tod{e}_1, \tod{e}_2}. 
\label{self20}
\end{equation}
\item $\toq{\Gamma}^0 \toq{L}^0$ is a positive operator from $\mathscr{H}$ endowed with the scalar product (\ref{scalar}) into itself:
\begin{equation}
\mmoy{\tod{e}, \toq{\Gamma}^0 \toq{L}^0 \tod{e}}  \geq 0, \quad \forall \; \tod{e} \in  \mathscr{H}.
\label{self3bis}
\end{equation}
\item 
$\toq{\Gamma}^0 \toq{L}^0$ is the orthogonal projector  from $\mathscr{H}$ onto $\mathscr{E}_0$ for the energetic scalar product (\ref{scalar}) 
\begin{equation}
\forall  \tod{e}  \in \mathscr{H}: \; \toq{\Gamma}^0 (\toq{L}^0\tod{e})\in \ \mathscr{E}_0, \; 
\toq{\Gamma}^0 \toq{L}^0 \tod{\eps}^* = \tod{\eps}^* \quad \forall \tod{\eps}^* \in \mathscr{E}_0, \;
\mmoy{\tod{e}- \toq{\Gamma}^0 \toq{L}^0 \tod{e}, \toq{\Gamma}^0 \toq{L}^0 \tod{e}}=0.
\label{prop5}
\end{equation}
As such it is a contraction on $\mathscr{H}$ and its  operator norm is exactly 1 when $\mathscr{H}$ is endowed with the energetic scalar product (\ref{scalar}).
\end{enumerate}
\section{Series expansions of Obnosov's analytical expression}
\label{dvpt}
The Obnosov relation (\ref{obnosov1}) 
\begin{equation}
  \eff{\toq{L}}(z) = \eff{z} \toq{I}, \quad  \eff{z} = \sqrt{   \frac{ 1 + 3 {z} }{ 3 + {z}}   },
  \nonumber
\end{equation}
can be expanded in power series with respect to the three variables of interest $t=z-1$, $\displaystyle t= Z= \frac{z-1}{z+1}$ or 
$\displaystyle  t=w= \frac{\sqrt{z}-1}{\sqrt{z}+1}$. 
\vskip 0.1cm
Two types of expansions are of interest:
\begin{equation}
\frac{\eff{z}}{z_0} = \sum_{k=0}^{\infty} b_k t^k, \quad \eff{z} = \sum_{k=0}^{\infty} d_k t^k.
\label{expansion}
\end{equation}
The first quantity will serve for comparison with the general expansion of ${\toq{L}^0}^{-1}\eff{\toq{L}}$, whereas the second expansion serves to assess the accuracy of the numerical simulations.
\paragraph{A useful trick}
First, because of the square root in (\ref{obnosov1}), it is noticed that an expansion of $(\eff{z}/z_0)^2$ or $\eff{z}^2$ will be easier to derive than that of  $\eff{z}/z_0$ or $\eff{z}$. Writing 
\begin{equation}
\left(\frac{\eff{z}}{z_0}\right)^2  = \sum_{n=0}^{\infty}{ a_n t^n}, \quad \frac{\eff{z}}{z_0}  = \sum_{k=0}^{\infty}{ b_k t^k}, \quad   \eff{z}^2  = \sum_{n=0}^{\infty}{ c_n t^n}, 
\quad \eff{z}  = \sum_{k=0}^{\infty}{ d_k t^k},
\label{an}
\end{equation}
the coefficients $(a_n)\vert_{n \in \mathbb{N}}$ and $(b_k)\vert_{k \in \mathbb{N}}$ are related by:
\begin{equation}
  a_n = \sum_{k=0}^{n}{b_k b_{n-k}},\; \text{i.e.} \; b_0^2 = a_0,\; 2 b_0 b_1 = a_1,\; 2 b_0 b_n + \sum_{k=1}^{n-1}{b_k b_{n-k}} = a_n \quad \forall n \ge 2.
\end{equation}
which gives
\begin{equation}
      b_0 = \sqrt{a_0},\;     b_1 = \frac{1}{2} \frac{a_1}{\sqrt{a_0}},\;      b_k = \frac{1}{2 b_0} ( a_k - \sum_{i=1}^{k-1}{b_i b_{k-i}} ) \quad \forall k \ge 2,
  \label{dev0}
\end{equation}
with similar relations between $c_n$ and $d_k$. 
\paragraph{Reference=matrix. $\displaystyle z_0=1,\; t={z-1}$}
\label{dvpt1}
With $t=z-1$, the relation (\ref{obnosov1}) can be rewritten as:
\begin{equation}\eff{z}^2= {   \frac{ 1 + 3 t/4 }{ 1 + t/4}   } 
 =  1   + \sum_{k=1}^{\infty}{ 2 (-1)^{k-1}\left(\frac{1}{4}\right)^k t^{k}} , 
 \label{dev1}
 \end{equation}
Therefore (with the notations of (\ref{an})), 
\begin{equation}
    a_0 = 1, \quad 
    a_k = 2 (-1)^{k-1} 4^{-k} \quad \forall k \ge 1 \\
\label{a1}
\end{equation}
and following (\ref{dev0}),
\begin{equation}
    b_0 = 1,\;\
    b_1 = \frac{1}{4},\;  
    b_k= (-1)^{k-1} \left( \frac{1}{4}\right)^{k} - \frac{1}{2} \sum_{i=1}^{k-1}{b_i b_{k-i}} \quad \forall n \ge 2\\
   \label{d1}
\end{equation}
The reference between the matrix ($z_0=1$), the coefficients $c_n$ and $d_k$ coincide with $a_n$ and $b_k$ respectively.
\paragraph{Reference=arithmetic mean. $\displaystyle z_0=\frac{z+1}{2},\; t=\frac{z-1}{z+1}$}  
\label{dvpt2}
Noting that 
\begin{equation}
  t =  \frac{z-1}{z+1}, \quad 
  z =  \frac{1+t}{1-t},\quad  z_0=\frac{z+1}{2} = \frac{1}{1-t}
\end{equation}
it follows that: 
\begin{equation}
  \left( \frac{\eff{z}}{z_0}\right)^2 = {   \frac{ 1 + t/2 }{ 1 - t/2}  (1-t)^2 } 
  =  1 -t-\frac{1}{2}t^2 + \sum_{k=3}^{\infty} 2^{1-k} t^{k}.
\label{dev2}  
\end{equation}
Therefore
\begin{equation}
    a_0 = 1, \; a_1=-1,\; a_2= -\frac{1}{2},\; a_k = 2^{1-k} \quad \forall k \ge 3.
\label{a2}        
\end{equation}
The $b_k$'s can then be obtained from (\ref{a2})  by (\ref{dev0}). 
Similarly \begin{equation}
  \eff{z}^2 = {   \frac{ 1 + t/2 }{ 1 - t/2}   } 
  =  1 +2 \sum_{k=1}^{\infty} \left(\frac{t}{2}\right)^{k}.
\label{dev2b}  
\end{equation}
Therefore
\begin{equation}
    c_0 = 1, \quad c_k = 2^{1-k} \quad \forall k \ge 1.
\label{c2}        
\end{equation}
 The $d_k$'s can then be obtained from (\ref{c2})  by (\ref{dev0}). 
\paragraph{Reference= geometric mean.  $\displaystyle z_0=\sqrt{z}, \; t=\frac{\sqrt{z}-1}{\sqrt{z}+1}$}  
\label{dvpt3}
With
\begin{equation}
  t =  \frac{\sqrt{z}-1}{\sqrt{z}+1},\quad  z =  \left( \frac{1+t}{1-t} \right) ^2 =z_0^2,
\end{equation}
 one can write
 \begin{equation}
  \left(\frac{\eff{z}}{z_0}\right)^2 =  \frac{ 1 + t + t^2 }{ 1 - t + t^2}  \left( \frac{1-t}{1+t} \right) ^2 = 
  \frac{ 1 -t^3 }{ 1 + t^3}   \frac{1-t}{1+t}.
  \end{equation}
Then, noting that 
$$  \frac{1-t}{1+t} = 1 + 2 \sum_{k=1}^{\infty} (-t)^{k},\quad \frac{1-t^3}{1+t^3} = 1 + 2 \sum_{k=1}^{\infty} (-t)^{3k}, $$
one arrives at
\begin{equation}
  \left(\frac{\eff{z}}{z_0}\right)^2   = \sum_{k=0}^{\infty}{ a_k t^k}, \quad a_0=a'_0a''_0,\quad   
  a_k = \sum_{i=0}^{k} a'_{i}a''_{k-i}, 
  \end{equation}
with 
\begin{equation}
\left.
\begin{array}{l}
a'_0=1,\quad  a'_i=2 (-1)^{i}\quad  \forall i \geq 1, \\[2ex]
 a''_0=1,\;  a''_1=0, \; a''_2=0, \;
a''_{3j}=2 (-1)^{3j},\; a''_{3j+1}=a''_{3j+2}= 0 \quad \forall j \geq 1
\end{array}
\right\}
\end{equation}
It is found that 
\begin{equation}
\left.
\begin{array}{l}
a_0=1,\;  a_1=-2 ,\; a_2=2,\\[2ex] 
a_{3j}=(-1)^{3j}(4j),\; a_{3j+1}=(-1)^{3j+1}(4j+2),\; a'_{3j+2}=  (-1)^{3j+2}(4j+2) \quad \forall j \geq 1.
\end{array}
\right\}
\end{equation}
Similarly
\begin{equation}
  \eff{z}^2 =  \frac{ 1 + t + t^2 }{ 1 - t + t^2}   = 
  \frac{ 1 -t^3 }{ 1 + t^3}   \frac{1+t}{1-t}.
  \end{equation}
\begin{equation}
\eff{z}^2   = \sum_{k=0}^{\infty}{ c_k t^k}, \quad c_0=c'_0c''_0,\quad   
  c_k = \sum_{i=0}^{k} c'_{i}c''_{k-i}, 
  \end{equation}
with 
\begin{equation}
\left.
\begin{array}{l}
c'_0=1,\quad  c'_i=2 \quad  \forall i \geq 1, \\[2ex]
 c''_0=1,\;  c''_1=0, \; c''_2=0, \;
c''_{3j}=2 (-1)^{3j},\; c''_{3j+1}=c''_{3j+2}= 0 \quad \forall j \geq 1
\end{array}
\right\}
\end{equation}
It is found that 
\begin{equation}
\left.
\begin{array}{l}
c_0=1,\;  c_1=2 ,\; c_2=2,\\[2ex] 
c_{3j}=0, \; c_{3j+1}=c_{3j+2}=2 (-1)^{3j} \quad \forall j \geq 1.
\end{array}
\right\}
\end{equation}
\newpage
\section{Comparison between the exact and numerical coefficients of the power series expansion for the Obnosov microstructure}
The theoretical and numerical coefficients $d_k$ for the expansion of $\eff{z}$ in powers of $\displaystyle t=\frac{z-1}{z+1}$ are given in the table below.   
\begin{center}
\begin{table}[!htpt]%
\centering
\begin{tabular*}{500pt}{@{\extracolsep\fill}lcccc@{\extracolsep\fill}}
\hline \\
    $k$ &$d_k$ (theoretical) &{${d}^{NS}_k$} (numerical) &$\displaystyle \frac{|d_k-{d^{NS}_k}|}{d_k}$ (\%)\\
\hline \\
 0  &1.000000000  	 &1.000000000           & 0.00000000      \\
1  &0.500000000	         &0.500000000           &0.00000000      \\
2  &0.125000000	         &0.125000000           &0.00000000      \\
3  &0.625000000E-01	 &0.625000000E-01       &0.00000000      \\
4  &0.234375000E-01	 &0.234375000E-01       &0.00000000      \\
5  &0.117187500E-01	 &0.117141463E-01       &0.03928491      \\
6  &0.488281250E-02	 &0.488051063E-02       &0.04714230      \\
7  &0.244140625E-02	 &0.243735662E-02       &0.16587284      \\
8  &0.106811523E-02	 &0.106609042E-02       &0.18956850      \\
9  &0.534057617E-03	 &0.535119341E-03       &0.19880327      \\
10 &0.240325928E-03	 &0.240856801E-03       &0.22089710      \\
11 &0.120162964E-03	 &0.125219714E-03       &4.20824340      \\
12 &0.550746918E-04	 &0.576030820E-04       &4.59083858      \\
13 &0.275373459E-04	 &0.342665252E-04       &24.4365573      \\
14 &0.127851963E-04	 &0.161497831E-04       &26.3162702      \\
15 &0.639259815E-05	 &0.133939187E-04       &109.522300      \\
16 &0.299653038E-05	 &0.649716274E-05       &116.822856      \\
17 &0.149826519E-05	 &0.822353999E-05       &448.870790      \\
18 &0.707514118E-06	 &0.407011143E-05       &475.269288      \\
19 &0.353757059E-06	 &0.666877547E-05       &1785.12859      \\
20 &0.168034603E-06	 &0.332550840E-05       &1879.06166      \\
21 &0.840173016E-07	 &0.600176714E-05       &7043.48953      \\
22 &0.400991667E-07	 &0.299895396E-05       &7378.84359      \\
23 &0.200495833E-07	 &0.559017674E-05       &27781.7602      \\
24 &0.960709201E-08	 &0.279466942E-05       &28989.6498      \\
25 &0.480354601E-08	 &0.527750045E-05       &109766.762      \\
\hline \\
  \end{tabular*}
\caption{Obnosov problem.  Discretization $128 \times 128$ pixels. MS-scheme. Comparison between the theoretical and numerical  coefficients $d_k$ of the power series with $\displaystyle t=\frac{z-1}{z+1}$. }
\label{table_dk}
\end{table}
\end{center}

\end{document}